\documentclass[12pt,francais]{smfart}
\usepackage{smfthm}
\usepackage[T1]{fontenc}
\usepackage[francais]{babel}
\usepackage[latin1]{inputenc}
\usepackage{amssymb,mathrsfs}
\let\mathcal\mathscr
\usepackage[matrix,arrow,curve,cmtip]{xy}
\CompileMatrices

 \calclayout

\makeatletter

 \numberwithin{paragraph}{subsection}
 
 \let\c@equation\c@paragraph
 \let\cl@equation\cl@paragraph
 \def\theequation {\thesubsection.\arabic{equation}}
 \let\subsection\Subsection
 \NumberTheoremsAs{paragraph}

 \def\th@plain{%
  \let\thm@indent\noindent
  \thm@headfont{\normalfont\sc}%
  \thm@notefont{\normalfont}
  \thm@preskip.5\linespacing
  \thm@postskip\thm@preskip
  \ifsmf@skippt
      \thm@headpunct{}\let\thmheadnl\@thmheadnl
      \global\smf@skipptfalse
  \else
       \thm@headpunct{\pointrait}
  \fi
  \itshape }

 \def\th@definition{%
  \th@plain
  \thm@headfont{\normalfont\itshape}\normalfont}

 \def\th@remark{%
   \th@plain
   \thm@headfont{\normalfont\itshape}\normalfont}
{\newcount\@hour\newcount\@minute
  \newcount\acl@temp\@hour=\time\@minute=\time
  \divide\@hour by 60\acl@temp=\@hour\multiply\acl@temp by60\relax
  \advance\@minute by-\acl@temp
  \global\edef\clocktime{\the\@hour:\ifnum\@minute<10 0\fi\the\@minute}}

\def\bibliosection{\@startsection{section}{1}%
  \z@{2\linespacing\@plus\linespacing}{.5\linespacing}%
  {\normalfont\bfseries\raggedright}}

\stepcounter{tocdepth}

     \def\l@part{\@tocline{-1}{0pt}{0pt}{}{\bfseries}}
     \def\l@section{\@tocline{1}{0pt}{0pt}{}{}}
     \def\l@subsection{\@tocline{2}{0pt}{12pt}{}{}}
     \def\l@subsubsection{\@tocline{3}{0pt}{}{}{}}
     \def\l@paragraph{\@tocline{5}{0pt}{}{}{}}

\def\baselinestretch{1.1}
\makeatother

\def\V{\mathbf{V}}
\def\P{\mathbf{P}}
\def\A{\mathbf{A}}
\def\R{\mathbf{R}}
\def\Q{\mathbf{Q}}

\def\C{\mathbf{C}}
\def\Z{\mathbf{Z}}

\let\bar\overline
\let\hat\widehat
\let\tilde\widetilde
\let\ra\rightarrow
\let\phi\varphi
\let\eps\varepsilon

\def\Aut{\operatorname{Aut}\nolimits}
\def\Sym{\operatorname{Sym}\nolimits}
\def\Spec{\operatorname{Spec}\nolimits}
\def\Proj{\operatorname{Proj}\nolimits}
\def\Qcoh{\operatorname{QCoh}\nolimits}
\def\Fib{\operatorname{Fib}\nolimits}
\def\Pic{\operatorname{Pic}\nolimits}
\def\div{\operatorname{div}\nolimits}
\def\rmH{{\mathrm H}}
\def\hatH{\widehat{\mathrm H}\vphantom{\mathrm H}}
\def\barH{\overline{\mathrm H}\vphantom{\mathrm H}}
\def\hatFib{\mathop{%
  \smash{\widehat{\operatorname{Fib}}}\vphantom{\operatorname{Fib}}}\nolimits}
\def\barPic{\mathop{%
  \smash{\overline{\operatorname{Pic}}}\vphantom{\operatorname{Pic}}}\nolimits}
\def\hatPic{\mathop{%
  \smash{\widehat{\operatorname{Pic}}}\vphantom{\operatorname{Pic}}}\nolimits}
\def\hCH{\widehat{\operatorname{CH}}\vphantom{\operatorname{CH}}}
\def\hdeg{\widehat{\operatorname{deg}}\,}
\def\vol{\operatorname{vol}}

\def\eff{\mathrm{eff}}

\def\Rep{\operatorname{Rep}}
\def\gm{\mathbf{G}_m}
\def\ga{\mathbf{G}_a}
\def\GL{\operatorname{GL}}
\def\Gal{\operatorname{Gal}}

\def\norm#1{\left\|{#1}\right\|}
\def\abs#1{\left|{#1}\right|}

\begin{document}

\title{Torseurs arithm\'etiques et espaces fibr\'es}

\author{Antoine Chambert-Loir}
\address{Institut de math\'ematiques de Jussieu\\ Boite 247 \\
4, place Jussieu \\ F-75252 Paris Cedex 05}
\email{chambert@math.jussieu.fr}

\author{Yuri Tschinkel}
\address{Department of Mathematics\\ U.I.C. \\ Chicago (IL) 60607 }
\email{yuri@math.uic.edu}

\date{Soumis sur l'archive \emph{alg-geom} le 4 janvier 1999}


\maketitle

{\def\baselinestretch{1}
\tableofcontents
}

\section*{Introduction}

Cet article est le premier d'une s\'erie 
dont le th\`eme principal est l'\'etude des hauteurs sur
certaines vari\'et\'es alg\'ebriques sur un corps de nombres. On
voudrait notamment comprendre la distribution des points rationnels
de hauteur born\'ee. 

Pr\'ecis\'ement, soient $X$ une vari\'et\'e alg\'ebrique projective lisse
sur un corps de nombres $F$, $\mathcal L$ un fibr\'e en droites
sur $X$ et $H_{\mathcal L}:X(\bar F)\ra\mathbf R_+^*$ une fonction
hauteur (exponentielle) pour $\mathcal L$. Si $U$ est un ouvert
de Zariski de $X$, on cherche \`a estimer le nombre
$$ N_U(\mathcal L,H) = \# \{ x\in U(F) \,;\, H_{\mathcal L}(x)\leq H \}$$
lorsque $H$ tend vers $+\infty$.
L'\'etude de nombreux exemples
a montr\'e que l'on peut s'attendre \`a un \'equivalent
de la forme
\begin{equation} \tag{$*$} \label{eq:maninpeyre} 
 N_U(\mathcal L,H) =
    \Theta(\mathcal L) H^{a(\mathcal L)} (\log H)^{b(\mathcal L)-1} 
    (1+o(1)), \quad H\ra+\infty 
\end{equation}
pour un ouvert $U$ convenable et lorsque par exemple $\mathcal L$ et
$\omega_X^{-1}$ (fibr\'e anticanonique) sont amples.
On a en effet un r\'esultat de ce genre lorsque 
$X$ est une vari\'et\'e de drapeaux~\cite{franke-m-t89},
une intersection compl\`ete lisse de bas degr\'e (m\'ethode du cercle),
une vari\'et\'e torique~\cite{batyrev-t98b},
une vari\'et\'e horosph\'erique~\cite{strauch-t98}, etc.
On dispose de plus d'une description conjecturale
assez pr\'ecise des constantes
$a(\mathcal L)$ et $b(\mathcal L)$ en termes du c\^one des diviseurs
effectifs~\cite{batyrev-m90}
ainsi que de la constante $\Theta(\mathcal L)$
(\cite{peyre95}, \cite{batyrev-t98}).

En fait, on \'etudie plut\^ot la \emph{fonction z\^eta des hauteurs},
d\'efinie par la s\'erie de Dirichlet
$$ Z_U(\mathcal L,s) =  \sum_{x\in U(F)} H_{\mathcal L}(x)^{-s} $$
\`a laquelle on applique des th\'eor\`emes taub\'eriens standard.
Sur cette s\'erie, on peut se poser les questions suivantes :
domaine de convergence, prolongement m\'eromorphe,
ordre du premier p\^ole, terme principal, sans oublier
la croissance dans les bandes verticales \`a gauche du premier p\^ole.
Cela permet de proposer des conjectures de pr\'ecision variable.

Il est naturel de vouloir tester la compatibilit\'e de cette conjecture
avec les constructions usuelles de la g\'eom\'etrie alg\'ebrique.
Par exemple, on n'arrive pas \`a d\'emontrer cette conjecture 
pour un \'eclatement $X'$ d'une vari\'et\'e $X$ pour laquelle 
cette conjecture est connue. M\^eme pour un \'eclatement
de 4 points dans le plan projectif, on n'a pas de r\'esultat complet !

\bigskip

Dans cet article, nous consid\'erons certaines fibrations
localement triviales construites de la fa\c con suivante.
Soient $G$ un groupe alg\'ebrique lin\'eaire sur $F$ agissant sur
une vari\'et\'e projective lisse $X$, $B$ une vari\'et\'e projective lisse
sur $F$ et $T$ un $G$-torseur sur $B$ localement trivial pour
la topologie de Zariski. Ces donn\'ees d\'efinissent une vari\'et\'e
alg\'ebrique projective $Y$ munie d'un morphisme $Y\ra B$ dont les
fibres sont isomorphes \`a $X$.
On donne au \S\,\ref{subsec:exemples} de nombreux exemples
{\og concrets\fg} de vari\'et\'es alg\'ebriques provenant d'une telle
construction.
Le c\oe ur du probl\`eme est de comprendre le comportement de
la fonction hauteur lorsqu'on passe d'une fibre \`a l'autre,
comportement vraiment non trivial bien qu'elles soient
toutes isomorphes.

Pour d\'efinir et \'etudier de fa\c con syst\'ematique
les fonctions hauteurs sur $Y$,
on est amen\'e \`a d\'egager de nouvelles notions dans
l'esprit de la g\'eom\'etrie d'Arakelov.
Apparaissent notamment les
notions de $G$-torseur arithm\'etique
au~\S\,\ref{defi:torsarith}, ainsi que 
la d\'efinition de la \emph{fonction $L$ d'Arakelov} attach\'ee
\`a un tel torseur arithm\'etique et \`a une fonction sur
le groupe ad\'elique $G(\mathbf A_F)$ invariante par $G(F)$ et
par un sous-groupe compact convenable (\S\,\ref{subsec:LArakelov}).
Elles g\'en\'eralisent les notions usuelles de fibr\'e inversible
m\'etris\'e ainsi que la fonction z\^eta des hauteurs introduits
par S.~Arakelov~\cite{arakelov74b}.

Ceci fait, on peut voir que les fonctions hauteurs d'une
fibre $Y_b$ de la projection $Y\ra B$ diff\`erent de la fonction
hauteur sur $X$ par
ce que nous appelons \emph{torsion ad\'elique}, dans laquelle
on retrouve explicitement la classe d'isomorphisme
du $G$-torseur arithm\'etique $T_b$ sur $F$ (\S\,\ref{subsec:torsion}).

Dans un deuxi\`eme article,
nous appliquerons ces consid\'erations g\'en\'erales au cas d'une
fibration en vari\'et\'es toriques provenant d'un torseur sous un
tore pour l'ouvert $U$ d\'efini par le tore.
Le principe de l'\'etude g\'en\'eralise~\cite{strauch-t98} 
et est le suivant.
On construit les hauteurs
\`a l'aide d'un prolongement du torseur g\'eom\'etrique
en un torseur arithm\'etique,
ce qui correspond en l'occurence
au choix de m\'etriques hermitiennes sur certains fibr\'es en droites.
On \'ecrit ensuite la fonction z\^eta
comme la somme des fonctions z\^eta des fibres 
$$ Z_U (\mathcal L,s)
= \sum_{b\in B(F)}  \sum_{x\in U_b(F)} H_{\mathcal L}(x)^{-s}
= \sum_{b\in B(F)} Z_{U_b}(\mathcal L|_{U_b},s). $$
Chaque $U_b$ est isomorphe au tore et on peut r\'ecrire la fonction z\^eta
des hauteurs de $U_b$ \`a l'aide de la formule de Poisson ad\'elique.
De cette fa\c con, la fonction z\^eta de $U$
appara\^{\i}t comme une int\'egrale
sur certains caract\`eres du tore ad\'elique
de la fonction $L$ d'Arakelov
d'un torseur arithm\'etique sur $B$.

Cette expression nous permettra d'\'etablir des th\'eor\`emes de mont\'ee
ou de descente : supposons que $B$ v\'erifie une conjecture,
alors $Y$  la v\'erifie ;
r\'eciproquement, supposons que $Y$ la v\'erifie,
alors, $B$ aussi. Bien s\^ur, la m\'ethode reprend les outils utilis\'es
dans la d\'emonstration de ces conjectures pour les vari\'et\'es toriques
(\cite{batyrev-t98b,batyrev-t95b,batyrev-t96}).

Alors que le pr\'esent article contient des consid\'erations
g\'en\'erales de {\og th\'eorie d'Arakelov \'equivariante\fg}
dont on peut esp\'erer qu'elles seront utiles dans d'autres contextes,
le deuxi\`eme verra intervenir des outils de th\'eorie analytique
des nombres (formule de Poisson, th\'eor\`eme des r\'esidus, 
estimations, etc.).

\bigskip
{\def\baselinestretch{1}\small
\def\thefootnote{\fnsymbol{footnote}}
\paragraph*{Remerciements}
Nous remercions J.-B.~Bost pour d'utiles discussions.
Pendant la pr\'eparation de cet article,
le second auteur\footnote[1]{partially supported by the N.S.A.}
\'etait invit\'e \`a l'I.H.E.S. 
et \`a Jussieu ; il est reconnaissant envers ces institutions
pour leur hospitalit\'e. \par }

\section*{Notations et conventions}

Si $\mathcal X$ est un sch\'ema, on d\'esigne par $\Qcoh(\mathcal X)$
et $\Fib_d(\mathcal X)$ les
cat\'egories des faisceaux quasi-coh\'erents (resp.\  des faisceaux
localement libres de rang $d$) sur $\mathcal X$.
On note $\Pic(\mathcal X)$ le groupe des classes d'isomorphisme
de faisceaux inversibles sur $\mathcal X$.
Si $\mathcal F$ est un faisceau quasi-coh\'erent sur $\mathcal X$,
on note $\V(\mathcal F)=\Spec\Sym \mathcal F$
et $\P(\mathcal F)=\Proj\Sym \mathcal F$ les fibr\'es vectoriels
et projectifs associ\'es \`a $\mathcal F$.

On note $\hatFib_d(\mathcal X)$ la cat\'egorie des fibr\'es vectoriels
hermitiens sur $\mathcal X$ (c'est-\`a-dire des faisceaux localement
libres de rang~$d$ munis d'une m\'etrique hermitienne
continue sur $\mathcal X(\C)$ et invariante par la conjugaison
complexe). On note $\hatPic(\mathcal X)$
le groupe des classes d'isomorphisme de fibr\'es en droites hermitiens
sur $\mathcal X$.

Si $\mathcal X$ est un $S$-sch\'ema, et si $\sigma\in S(\C)$,
on d\'esigne par $\mathcal X_\sigma$
le $\C$-sch\'ema $\mathcal X\times_\sigma \C$.
Cette notation servira lorsque $S$ est le spectre d'un localis\'e
de l'anneau des entiers d'un corps de nombres $F$, de sorte
que $\sigma$ n'est autre qu'un plongement de $F$ dans $\C$.

Si $G$ est un sch\'ema en groupes sur $S$,
$X^*(G)$ d\'esigne le groupe des $S$-homomorphismes
$G\ra\gm$ (caract\`eres alg\'ebriques).

Si $\mathcal X/S$ est lisse, le faisceau canonique de $\mathcal X/S$,
not\'e $\omega_{\mathcal X/S}$
est la puissance ext\'erieure maximale de $\Omega^1_{\mathcal X/S}$.

\section{Torseurs arithm\'etiques}

\subsection{D\'efinitions}

Rappelons la d\'efinition d'un torseur en g\'eom\'etrie alg\'ebrique.
\begin{defi}\label{defi:torseur}
Soient $S$ un sch\'ema, $\mathcal B$ un $S$-sch\'ema et
$G$ un $S$-sch\'ema en groupes plat et localement
de pr\'esentation finie.

Un \emph{$G$-torseur} sur un $\mathcal B$
est un $\mathcal B$-sch\'ema $\pi:\mathcal T\ra \mathcal B$
fid\`element plat et localement de pr\'esentation finie
muni d'une action de $G$ au-dessus de $\mathcal B$,
$m:G\times_S \mathcal T \ra \mathcal T$, de sorte que 
le morphisme 
$$
(m,p_2):G\times_S \mathcal T \ra \mathcal T\times_{\mathcal B} \mathcal T
$$
soit un isomorphisme.
On le suppose de plus \emph{localement trivial pour la topologie de
Zariski.}

On note $\rmH^1(\mathcal B,G)$ l'ensemble
des classes d'isomorphisme de $G$-torseurs sur $\mathcal B$.
\end{defi}

\begin{enonce}{Situation} \label{situ:GK}
Supposons que $S$ est le spectre de l'anneau des entiers d'un corps de nombres
$F$ et que $G$ est un $S$-sch\'ema en groupes lin\'eaire connexe plat
et de pr\'esentation finie.
Fixons pour tout plongement complexe de $F$, $\sigma\in S(\C)$, un sous-groupe
compact maximal $K_\sigma$ de $G(\C)$ et notons
$\mathbf K_\infty$ la collection $(K_\sigma)_{\sigma}$.
On suppose que pour deux plongements complexes conjugu\'es, les sous-groupes
compacts maximaux correspondants sont \'echang\'es par la conjugaison complexe.
\end{enonce}

\begin{defi}\label{defi:torsarith}
On appelle \emph{$(G,\mathbf K_\infty)$-torseur arithm\'etique} sur $\mathcal B$
la donn\'ee d'un $G$-torseur $\mathcal T$ sur $\mathcal B$
ainsi que pour tout $\sigma\in S(\C)$,
d'une section du $K_\sigma\backslash G_\sigma (\C)$-fibr\'e 
sur $\mathcal B_\sigma(\C)$
quotient \`a $\mathcal T_\sigma(\C)$ par l'action de $K_\sigma$.
On suppose de plus que pour deux plongements complexes conjugu\'es,
les sections sont \'echang\'ees par la conjugaison complexe.

On note
$\hatH^1(\mathcal B,(G,\mathbf K_\infty))$ l'ensemble
des classes d'isomorphisme de $(G,\mathbf K_\infty)$-torseurs arithm\'etiques
sur $\mathcal B$.

On note aussi 
$\hatH^0 (\mathcal B,(G,\mathbf K_\infty))$ l'ensemble des
sections $g\in \rmH^0(\mathcal B,G)$ telles que pour toute place \`a l'infini
$\sigma$,
$g$ d\'efinisse une section $\mathcal B_\sigma(\C)\ra \mathbf K_\sigma$.
\end{defi}

\begin{rema}
Se donner une section du $K_\sigma\backslash G_\sigma(\C)$-fibr\'e
associ\'e \`a $\mathcal T_\sigma(\C)$ sur $\mathcal B_\sigma(\C)$
revient \`a fixer dans un recouvrement ouvert $(U_i)$ pour la topologie
complexe les fonctions de transition $g_{ij}\in\Gamma(U_i\cap U_j,G)$
\`a valeurs dans $K_\sigma$. Il en existe car $G_\sigma(\C)$ est
hom\'eomorphe au produit de $K_\sigma$ par un $\R$-espace vectoriel
de dimension finie, cf.\ par exemple~\cite{borel50}.

D'autre part, on choisit dans cet article de supposer la section
continue. Dans certaines situations, il pourrait \^etre judicieux
de la supposer ind\'efiniment diff\'erentiable.
\end{rema}

La d\'ependance de cette notion en les sous-groupes maximaux fix\'es
est la suivante : 
toute famille $(x_\sigma)\in\prod_\sigma G_\sigma(\C)$ telle
que $K'_\sigma=x_\sigma K_\sigma x_\sigma^{-1}$ d\'etermine
une bijection canonique 
$$\hatH^1 (\mathcal B,(G,\mathbf K_\infty))
\simeq \hatH^1 (\mathcal B,(G,\mathbf K_\infty')). $$
(Rappelons que deux sous-groupes compacts maximaux sont conjugu\'es.)

\paragraph{Variante ad\'elique}

Il existe une variante ad\'elique des consid\'erations pr\'ec\'edentes qui supprime
en apparence la r\'ef\'erence \`a un mod\`ele sur $\Spec\mathfrak o_F$.
En effet, si $\mathcal B$ est propre sur $\Spec\mathfrak o_F$,
remarquons que pour toute place finie de $F$, un 
$G$-torseur arithm\'etique sur $\mathcal B$ induit
une section du morphisme
$G(\mathfrak o_v)\backslash \mathcal T(F_v)\ra \mathcal B(F_v)$.

\begin{defi}
Soit $G_F$ un $F$-sch\'ema en groupes de type fini
et fixons un sous-groupe compact maximal\footnote{%
Cela signifiera pour nous que les $K_v$ sont des sous-groupes
compacts ouverts aux places finies, et maximaux aux places infinies.}
$\mathbf K=\prod_v K_v$ 
du groupe ad\'elique $G(\A_F)$.
Soit $\mathcal B_F$ un $F$-sch\'ema propre.

On appelle \emph{$(G_F,\mathbf K)$-torseur ad\'elique}
sur $\mathcal B_F$ la donn\'ee d'un $G_F$-torseur
$\mathcal T_F\ra\mathcal B_F$,
ainsi que pour toute place $v$ de $F$, d'une section
continue de $K_v\backslash \mathcal T_F(F_v) \ra \mathcal B_F(F_v)$.
On suppose de plus qu'il existe un ouvert non vide $U$
de $\Spec\mathfrak o_F$,
un $U$-sch\'ema en groupes plat et de pr\'esentation fini $G$,
un $U$-sch\'ema $\mathcal B$ propre, plat et de type fini,
ainsi qu'un $G$-torseur $\mathcal T\ra\mathcal B$  
qui prolongent respectivement $G_F$, $\mathcal B_F$ et $\mathcal T_F$
et v\'erifiant :
pour toute place finie  $v$ de $F$ dominant $U$,
$G(\mathfrak o_v)=K_v$ et la section continue de
$K_v \backslash \mathcal T_F(F_v)\ra\mathcal B_F(F_v)$
est celle fournie par le mod\`ele
$\mathcal T\ra\mathcal B$.

On note $\barH^1(\mathcal B_F,(G_F,\mathbf K))$ l'ensemble des classes
d'isomorphisme de $(G_F,\mathbf K)$-torseurs ad\'eliques sur $\mathcal B_F$.
\end{defi}

Bien s\^ur, si $\mathcal B$ est un $\mathfrak o_F$-sch\'ema propre
et $G$ un $\mathfrak o_F$-sch\'ema en groupes plat et de pr\'esentation
finie, tout $(G,\mathbf K_\infty)$-torseur arithm\'etique sur $\mathcal B$
d\'efinit un $(G_F,\mathbf K)$-torseur ad\'elique o\`u
$\mathbf K $ est le compact ad\'elique
$\prod_{\text{$v$ finie}} G(\mathfrak o_v) \prod_\sigma K_\sigma$.

\paragraph{Exemples}
a) Quand $G=\GL(d)$,
le torseur $\mathcal T$ correspond naturellement
\`a la donn\'ee d'un fibr\'e vectoriel $\mathcal E$ de rang $d$ sur $\mathcal B$
par la formule
$\mathcal T = \mathbf{Isom}(\mathcal O_{\mathcal B}^n,\mathcal E)$.
Si l'on choisit $K_\sigma=\mathrm U(d)$,
une section du $\mathrm{U}(d,\C)\backslash \GL(d,\C)$-fibr\'e associ\'e
correspond \`a une m\'etrique hermitienne (continue) sur $\mathcal E$.
Ainsi, les $(\GL(d),\mathrm U(d))$-torseurs arithm\'etiques sont en bijection
naturelle avec les fibr\'es vectoriels hermitiens.

b) En particulier, lorsque $G=\gm$, la famille des sous-groupes
compacts maximaux $\mathbf K_\infty$ est canoniquement d\'efinie 
(ce qui permet de les omettre dans la notation) et
$\hatH^1(\mathcal B,\gm)=\hatPic(\mathcal B)$,
le groupe des classes d'isomorphisme
de fibr\'es en droites sur $\mathcal B$ munis d'une m\'etrique hermitienne
continue compatible \`a la conjugaison complexe. 
Les $\gm$-torseurs ad\'eliques s'identifient de m\^eme aux fibr\'es
inversibles munis d'une m\'etrique ad\'elique. Nous rappelons
cette th\'eorie au paragraphe~\ref{subsec:metriqueadelique}

c) Dans ce texte, nous ne consid\'erons que des $G$-torseurs localement
triviaux pour la topologie de Zariski. 
N\'eanmoins, lorsque $G/S$ est un $S$-sch\'ema ab\'elien, un exemple
de $G$-torseur localement trivial pour la topologie \'etale sur $\mathcal B$
est fourni par un sch\'ema ab\'elien $\mathcal A/\mathcal B$
obtenu par torsion de $G/S$, c'est-\`a-dire tel qu'il existe
un rev\^etement \'etale $\mathcal B'\ra\mathcal B$ de sorte
que $\mathcal A\times_{\mathcal B}\mathcal B'$ soit isomorphe
\`a $G\times_S \mathcal B'$
(famille de sch\'emas ab\'eliens \`a module constant).
De tels exemples devraient bien s\^ur faire partie d'une \'etude 
plus g\'en\'erale de la g\'eom\'etrie d'Arakelov des torseurs que
nous reportons \`a une occasion ult\'erieure.

\subsection{Propri\'et\'es}

Les ensembles de classes d'isomorphisme de $(G,\mathbf K_\infty)$-torseurs
arithm\'etiques v\'erifient un certain nombre  de propri\'et\'es formelles,
dont les analogues alg\'ebriques sont bien connus.
Leur d\'emonstration est standard et laiss\'ee au lecteur.

\begin{prop} \label{prop:H/hatH}
L'oubli de la structure arithm\'etique induit une application
$$\hatH^1 (\mathcal B,(G,\mathbf K_\infty)) \ra \rmH^1 (\mathcal B,G). $$
On a aussi une suite exacte d'ensembles point\'es :
\begin{multline*}
1 \ra \hatH^0(\mathcal B, (G,\mathbf K_\infty))
     \ra \rmH^0(\mathcal B, G)
    \ra \left(\bigoplus_\sigma \Gamma(\mathcal B_\sigma(\C),
K_\sigma\backslash G_\sigma(\C))\right)^{F_\infty} \ra \\
    \ra \hatH^1(\mathcal B,(G,\mathbf K_\infty))
    \ra \rmH^1(\mathcal B,G) \ra 1. 
\end{multline*}
\end{prop}
\begin{rema}
Lorsque $G=\gm$, en identifiant $\gm(\C)/K$ \`a $\R_+^*$,
nous retrouvons
la suite exacte bien connue pour
$\hatPic$ et $\Pic$ (cf.~\cite{gillet-s90}, 3.3.5 ou 3.4.2).

D'autre part, on devrait pouvoir interpr\'eter cette suite exacte
\`a l'aide de la  \emph{mapping cylinder category} introduite par
S.~Lichtenbaum dans son \'etude des valeurs sp\'eciales
des fonctions z\^eta des corps de nombres.
En effet, cette cat\'egorie est (?!)
la cat\'egorie des faisceaux en groupes ab\'eliens
sur, disons $\Spec\Z\cup\{\infty\}$.
\end{rema}

\begin{prop}
Supposons que le groupe $G$ est commutatif. Alors,
les sous-groupes compacts maximaux
sont uniques et l'ensemble $\hatH^1 (\mathcal B,(G,\mathbf K_\infty))$
h\'erite d'une structure de groupe ab\'elien
compatible avec la structure de groupe ab\'elien sur
$\rmH^1 (\mathcal B,G)$.
Dans ce cas, la suite exacte~\ref{prop:H/hatH} est une suite exacte
de groupes ab\'eliens.
\end{prop}

\begin{prop}
\emph{(Changement de base)}
Tout morphisme de $S$-sch\'emas $\mathcal B'\ra\mathcal B$ induit
un foncteur des $(G,\mathbf K_\infty)$-torseurs arithm\'etiques sur $\mathcal B$
vers les $(G,\mathbf K_\infty)$-torseurs arithm\'etiques sur $\mathcal B'$,
compatible \`a l'oubli des structures arithm\'etiques et aux
classes d'isomorphisme.

\emph{(Changement du corps de base)}
Si $F'$ est une extension de $F$, $S'=\Spec\mathfrak o_{F'}$
et si on choisit pour tout plongement complexe $\sigma'$ de $F'$
$K_{\sigma'}=K_{\sigma'|_F}$, on dispose d'un foncteur
des $(G,\mathbf K_\infty)$-torseurs arithm\'etiques sur $\mathcal B$
vers les $(G\times_S S',\mathbf K_\infty)$-torseurs arithm\'etiques
sur $\mathcal B\times_S S'$, compatible \`a l'oubli
des structures arithm\'etiques et aux classes d'isomorphisme.
\end{prop}

\begin{prop}
\emph{(Changement de groupe)}
Si $p: G\ra G'$ est un morphisme de $S$-sch\'emas en groupes
et que les sous-groupes compacts maximaux $\mathbf K_\infty$
et $\mathbf K'_\infty$ sont choisis de sorte que
pour tout plongement complexe $\sigma$,
tel que $p(K'_{\sigma})\subset K'_{\sigma}$,
il y a un foncteur des $(G,\mathbf K_{\infty})$-torseurs arithm\'etiques
vers les $(G',\mathbf K'_{\infty})$-torseurs arithm\'etiques, compatible
\`a l'oubli des structures arithm\'etiques et aux classes d'isomorphisme.

\emph{(Suite exacte courte)}
Soit $$1\ra G''\stackrel\iota\ra G\stackrel p\ra G'\ra 1$$
une suite exacte de $S$-sch\'emas en groupes.
Soient $\mathbf K_\infty$, $\mathbf K_\infty'$ et $\mathbf K_\infty''$ des familles de sous-groupes
compacts maximaux pour $G$, $G'$ et $G''$ aux places
archim\'ediennes choisis de sorte que $K''_\sigma=\iota^{-1}(K_\sigma)$
et $p(K_\sigma) = K'_\sigma$ pour toute place $\sigma$.

Si $p$ admet localement une section (comme $S$-sch\'ema), alors on a une suite
exacte courte canonique d'ensembles point\'es :
\begin{multline*}
 1 \ra
{\hatH^0(\mathcal B,(G'',\mathbf K_\infty'')) } \stackrel{\iota}\ra 
{\hatH^0(\mathcal B,(G,\mathbf K_\infty)) } \stackrel{p} \ra 
{\hatH^0(\mathcal B,(G',\mathbf K_\infty')) } \stackrel{\delta} \ra\\
\ra
{\hatH^1(\mathcal B,(G'',\mathbf K_\infty'')) } \stackrel{\iota}\ra
{\hatH^1(\mathcal B,(G,\mathbf K_\infty)) } \stackrel{\pi}\ra
{\hatH^1(\mathcal B,(G',\mathbf K_\infty')). }
\end{multline*}
\end{prop}

\medskip

Sur $\Spec\mathfrak o_F$, l'ensemble des classes d'isomorphisme
de $(G,\mathbf K_\infty)$-torseurs arithm\'etiques a une description
tr\`es simple, similaire \`a la description classique des
classes d'isomorphisme de
$G$-torseurs sur une courbe projective sur un corps fini.
Cela g\'en\'eralise la description analogue
du groupe $\hatPic(\Spec\mathfrak o_F)$
(cf.~\cite{gillet-s90}, 3.4.3, p.~131, o\`u le groupe correspondant est
not\'e $\hCH^1(\Spec\mathfrak o_F)$).
\begin{prop} \label{prop:H1-adelique}
On a des isomorphismes canoniques
$$ \hatH^1 (\Spec\mathfrak o_F, (G,\mathbf K_\infty)) 
   \simeq G(F) \backslash G(\A_F) / \mathbf K_G, $$
o\`u $\mathbf K_G$ d\'esigne le produit
$\prod\limits_{\text{$v$ finie}} G(\mathfrak o_v)
\prod\limits_{\text{$\sigma $ infinie}} K_\sigma $.

De m\^eme, pour un sous-groupe compact maximal $\mathbf K$ de $G(\A_F)$,
on a un isomorphisme canonique
$$ \barH^1(\Spec F, (G_F,\mathbf K))
\simeq G(F)\backslash G(\A_F) / \mathbf K . $$
\end{prop}
\begin{proof}
Soit $\widehat{\mathcal T}$ un $(G,\mathbf K)$-torseur arithm\'etique sur
$\Spec(\mathfrak o_F)$, localement trivial pour la topologie de Zariski.
Commen\c cons par fixer un section $\tau_F\in\mathcal T(F)$.
Si $v$ est une place finie de $F$, comme $\rmH^1 (\Spec\mathfrak o_v, G)=0$,
il existe une section $\tau_v\in\mathcal T(\mathfrak o_v)$, 
unique modulo l'action de $G(\mathfrak o_v)$. Cette section
se relie \`a $\tau_F$ par un \'element bien d\'efini
$g_v\in G(F_v)/G(\mathfrak o_v)$ tel que
$g_v^{-1}\cdot \tau_F= \tau_v$.
Comme $\tau_F$ s'\'etend en une section de $\mathcal T$ sur un ouvert
de $\Spec\mathfrak o_F$, on a $g_v\in\mathbf K_v$ pour presque toute place $v$.
D'autre part, si $\sigma$ est une place infinie, la section
de $K_\sigma\backslash \mathcal T(\C)$ donn\'ee par 
la structure de $(G,\mathbf K_\infty)$-torseur arithm\'etique
est de la forme $K_\sigma g_\sigma ^{-1} \tau_F$,
pour un unique $g_\sigma\in G(\C)/K_\sigma$.
On a ainsi d\'efini un \'element $\mathbf g$ dans $G(\A_F)/\mathbf K_G$.
Il d\'epend de la section $\tau_F$, mais si on choisit une autre
section, elle sera de la forme $g_F \tau_F$, ce qui revient \`a changer
l'\'el\'ement $\mathbf g$ par $g_F^{-1} \mathbf g$.
Nous avons donc attach\'e au $(G,\mathbf K_\infty)$-torseur arithm\'etique
$\widehat{\mathcal T}$
un \'el\'ement dans $G(F)\backslash G(\A_F)/\mathbf K_G$
qui visiblement ne d\'epend que de la classe d'isomorphisme de
$\widehat{\mathcal T}$.

Pour la bijection r\'eciproque, on choisit un repr\'esentant
de $\mathbf g\in G(F)\backslash G(\A_F)/\mathbf K_G$
o\`u pour toute place finie $v$, $g_v\in G(F)$,
et o\`u presque tous les $g_v$ valent $1$.
Soit alors $U$ le plus grand ouvert de $\Spec\mathfrak o_F$
tel que pour toute place finie $v$, $g_v\in G(U)$ ;
si $v$ est une place finie qui ne domine pas $U$,
soit $U_v=U\cup\{v\}$.
On d\'efinit un $G$-torseur $\mathcal T$ sur $\Spec\mathfrak o_F$
comme isomorphe \`a $G$ sur $U$ et sur chaque $U_v$, les isomorphismes
de transition \'etant fix\'es par l'isomorphisme entre
$\mathcal T|_U=G|_U$ et $\mathcal T|_{U_v}\times U=G|_U$ induit
par la multiplication \`a gauche par $g_{v}^{-1}$.
On munit ce $G$-torseur de la $K_\sigma$-classe \`a gauche
$K_\sigma g_\sigma^{-1}$ dans la trivialisation 
canonique sur l'ouvert $U$ qui contient $\Spec F$,
d'o\`u un $(G,\mathbf K_\infty)$-torseur arithm\'etique
sur $\Spec\mathfrak o_F$.

On laisse au lecteur le soin de v\'erifier plus en d\'etail
que la classe d'isomorphisme du $(G,\mathbf K_\infty)$-torseur
arithm\'etique ainsi construit est ind\'ependante du repr\'esentant
choisi, et que cela d\'efinit effectivement la bijection r\'eciproque voulue.

La variante ad\'elique $\barH^1(\Spec F, (G_F,\mathbf K))$
se traite de m\^eme (et plus facilement car on n'a pas de
torseur \`a construire !).
\end{proof}
\begin{rema}
On aurait aussi pu construire le $G$-torseur $\mathcal T$ associ\'e \`a
un point ad\'elique $(g_v)$ en d\'ecr\'etant que les sections
de $\mathcal T$ sur un ouvert $U$ de $\Spec \mathfrak o_F$
sont les $\gamma\in G(F)$ tels que pour toute place finie
$v$ dominant $U$, $\gamma g_v\in G(\mathfrak o_v)$.
\end{rema}

\subsection{M\'etriques ad\'eliques}
\label{subsec:metriqueadelique}

Pour la commodit\'e du lecteur, nous rappelons la th\'eorie
des m\'etriques ad\'eliques sur les fibr\'es en droites. C'est
un cas particulier bien connu des constructions pr\'ec\'edentes
lorsque le groupe est $\gm$,
mais l'exposer nous permettra de fixer quelques notations.

\begin{defi}
Soient $F$ un corps valu\'e, $X$ un sch\'ema de type fini sur $F$
et $\mathcal L$ un fibr\'e en droites sur $X$.
Une m\'etrique sur $\mathcal L$ est une application continue
$\V(\mathcal L^\vee)(F)\ra \R_+$ de sorte que pour tout $x\in X(F)$,
la restriction de cette application \`a la fibre en $x$
(identifi\'ee naturellement \`a $F$) soit une norme.
\end{defi}

Soient $F$ un corps de nombres, $X$ un sch\'ema projectif sur $F$
et $\mathcal L$ un fibr\'e en droites sur $X$.
La donn\'ee d'un sch\'ema projectif et plat $\tilde X$ sur le spectre
$S=\Spec \mathfrak o_F$
de l'anneau des entiers de $F$ dont la fibre g\'en\'erique est $X$
d\'efinit pour toute place non-archim\'edienne
$v$ de $F$ une m\'etrique sur le fibr\'e en droites $\mathcal L\otimes F_v$
sur $X\times F_v$.

\begin{defi}
On appelle
\emph{m\'etrique ad\'elique} sur $\mathcal L$
toute collection de m\'etriques $(\norm{\cdot}_v)_v$ sur
$\mathcal L\otimes F_v$ pour toutes les places $v$ de $F$
qui est obtenue de cette fa\c con pour presque toutes les places
(non-archim\'ediennes) de $F$.

On note $\barPic(X)=\barH^1(X,\gm)$ le groupe des classes d'isomorphisme de
fibr\'es en droites sur $X$ munis de m\'etriques ad\'eliques.
\end{defi}

Donnons nous une m\'etrique ad\'elique sur $\mathcal L$.
Tout morphisme $f:Y\ra X$ de $F$-sch\'emas projectifs fournit
par image r\'eciproque
une m\'etrique ad\'elique sur $f^*\mathcal L$.
Si $Y$ n'est pas projective, on obtient tout de m\^eme de la sorte
une collection de m\'etriques pour toutes les places de $F$.

\begin{defi} \label{hauteuradelique}
Si $\overline{\mathcal L}=(\mathcal L,(\norm{\cdot}_v)_v)$
est un fibr\'e en droites sur $X$ muni d'une m\'etrique ad\'elique,
on appelle \emph{fonction hauteur} (exponentielle)
associ\'ee \`a $\overline{\mathcal L}$ la fonction
$$ H(\overline{\mathcal L}; \cdot) : X(F) \ra \R_+ ,
\quad x\mapsto \prod_{v} \norm{\mathsf s}_v (x)^{-1} , $$
$\mathsf s$ \'etant une section non nulle arbitraire
de $\mathcal L|_x\simeq F$.

Si $\mathsf s$ est une section globale non nulle de $\mathcal L$,
on d\'efinit une
\emph{fonction hauteur} (exponentielle)
\emph{sur les points ad\'eliques de $X$} en posant
$$ H(\overline{\mathcal L},\mathsf s;\cdot) :
    X(\A_F)\setminus |\div(s)|\ra\R_+,
   \quad  \mathbf x=(x_v)_v\mapsto \prod_v \norm{\mathsf s}_v (x_v)^{-1}.
$$
\end{defi}
\noindent
(Dans les deux cas, le produit converge en effet car il n'y a qu'un
nombre fini de termes diff\'erents de $1$.)
D'autre part,
elle est multiplicative en le fibr\'e en droites (resp.\ en la section),
ce qui permettra de l'\'etendre aux groupes de Picard tensoris\'es par $\C$.

Comme on a un isomorphisme
canonique $\barPic(\Spec F)=\hatPic(\Spec\mathfrak o_F)$,
on remarque que
$$ H(\overline{\mathcal L};x) = \exp(\hdeg \overline{\mathcal L}|_x) $$
o\`u $\hdeg : \hatPic(\Spec\mathfrak o_F)\ra\R$ est l'homomorphisme
{\og degr\'e arithm\'etique\fg}
d\'efini dans \cite{gillet-s90}, 3.4.3, p.~131.
Par l'isomorphisme de \emph{loc.\ cit.,} 
$$ \hatPic(\Spec\mathfrak o_F)\ra F^\times\backslash \A_F^\times / K , $$
$\exp\circ \hdeg$ correspond \`a l'inverse de la norme.

\begin{defi}
Soit $X$ une vari\'et\'e sur $F$, $\bar{\mathcal L}\in\barPic(X)_\C$
(le groupe des fibr\'es inversibles sur $X$ munis d'une m\'etrique ad\'elique
tensoris\'e par $\C$).
Si $U\subset X$ est un ouvert de Zariski, on appelle
fonction z\^eta des hauteurs de $U$ en $\bar{\mathcal L}$ la somme
$$ Z_U(\bar{\mathcal L}) = \sum_{x\in U(F)} H(\bar{\mathcal L};x)^{-1} $$
quand elle existe.
\end{defi}

\begin{rema}
La convergence absolue de la s\'erie ne d\'epend que de la partie
r\'eelle de $\bar{\mathcal L}$ dans $\Pic(X)_\R$ (on peut comparer deux
m\'etriques ad\'eliques). De plus, l'ensemble des $\mathcal L\in\Pic (X)_\R$
pour lesquels la s\'erie converge est une partie convexe
(in\'egalit\'e arithm\'etico-g\'eom\'etrique).
Enfin, si $\mathcal L$ est ample, alors $Z_U(s\bar{\mathcal L})$ converge
pour $\Re(s)$ assez grand
et d\'efinit une fonction analytique de $s$, not\'ee $Z_U(\bar{\mathcal L},s)$
dans l'introduction.
\end{rema}

Les consid\'erations analogues sont \'evidemment valables pour le groupe
de Picard--Arakelov $\hatPic(\mathcal X)$ d'un mod\`ele propre et plat
$\mathcal X$ de $X$ sur $\Spec\mathfrak o_F$.

\begin{enonce}[remark]{Exemple}
Lorsque $X$ est une vari\'et\'e torique, $\Pic(X)_\R$ est un espace
vectoriel de dimension finie et
il y a des m\'etriques canoniques
sur les fibr\'es en droites sur $X$ (cf.~\cite{batyrev-t95b}),
d'o\`u un homomorphisme canonique $\Pic(X)_\C\ra\barPic(X)_\C$.
Batyrev et Tschinkel ont montr\'e dans~\cite{batyrev-t98b}
que la s\'erie d\'efinissant la fonction z\^eta des hauteurs du tore
converge d\`es
que $\mathcal L\otimes \omega_{X}$ est dans l'int\'erieur
du c\^one effectif $\Lambda_{\eff}^\circ(X)\subset
\Pic(X)_\R$, le fibr\'e en droites $\mathcal L$ \'etant muni
de sa m\'etrique ad\'elique canonique. Elle d\'efinit m\^eme une fonction
holomorphe dans le tube sur ce c\^one.
\end{enonce}


\subsection{Fonctions $L$ d'Arakelov}
\label{subsec:LArakelov}

On se place dans la situation~\ref{situ:GK}.
Soient $\mathcal B$ un sch\'ema propre et fid\`element plat
sur $S=\Spec\mathfrak o_F$
et $\hat{\mathcal T}$ un $(G,\mathbf K_\infty)$-torseur arithm\'etique sur
$\mathcal B$.

Pour tout $b\in\mathcal B(F)$, il existe une unique
section $\eps_b:\Spec\mathfrak o_F\ra\mathcal B$ qui prolonge $b$.
On dispose ainsi d'un $(G,\mathbf K_\infty)$-torseur arithm\'etique
$\eps_b^*\hat{\mathcal T}$ sur $\Spec\mathfrak o_F$
que l'on notera $\hat{\mathcal T}|_b$.
En particulier, si $\Phi$ est une fonction \`a valeurs
complexes sur
$$
 G(F)\backslash G(\A_F) / \mathbf K_G
\simeq 
\hatH^1 (\Spec\mathfrak o_F, (G,\mathbf K_\infty)), $$
la composition
$$ 
\hatH^1 (\mathcal B, (G,\mathbf K_\infty))
\xrightarrow{\eps_b}
\hatH^1 (\Spec\mathfrak o_F, (G,\mathbf K_\infty))
\xrightarrow{\chi}
\C
$$
d\'efinit un nombre complexe $\Phi(\hat{\mathcal T}|_b)$.

\begin{defi} \label{defi:LArakelov}
Soient $\Phi$ une fonction sur $G(F)\backslash G(\A_F)/\mathbf K_G$
et $U$ une partie de $\mathcal B(F)$.
On appelle \emph{fonction $L$ d'Arakelov} l'expression
$$ L(\hat{\mathcal T},U,\chi) = \sum_{b\in U\subset\mathcal B(F)}
\Phi(\hat{\mathcal T}|_b), $$
quand la s\'erie converge (absolument).
\end{defi}

\paragraph{Exemple}
Soit $\hat{\mathcal L}\in\hatPic(\mathcal B)$ identifi\'e au
$\gm$-torseur arithm\'etique qu'il d\'efinit.
Si $U$ est l'ensemble des points rationnels d'un ouvert de $\mathcal B$,
la fonction $L$ d'Arakelov $L(\hat{\mathcal L},U,\norm{\cdot}^s)$ 
d\'efinie au \S\,\ref{subsec:LArakelov} ($\norm{\cdot}$
d\'esigne la norme ad\'elique) n'est autre que la fonction z\^eta
d'Arakelov $Z_U(\mathcal L,s)$,
introduite par Arakelov et largement \'etudi\'ee depuis.

En revanche, lorsque $\chi$ est un quasi-caract\`ere arbitraire de
$\hatPic(\Spec\mathfrak o_F)$ (pour la topologie ad\'elique),
on obtient un nouvel invariant $L(\hat{\mathcal L},U,\chi)$
dont l'importance appara\^{\i}tra \`a la fin de cet article.

\begin{rema}
Bien entendu, on d\'efinit de la m\^eme fa\c con
une fonction $L$ d'Arakelov,
$L(\bar{\mathcal T},U,\Phi)$
attach\'ee \`a un torseur ad\'elique $\bar{\mathcal T}$ sur $\mathcal B$
(sur $F$) et \`a une fonction $\Phi$ sur $G(F)\backslash G(\A_F)/\mathbf K$.
\end{rema}

\paragraph{Fonctions $\theta$ et $\zeta$}
Dans la suite de cette section, on suppose pour simplifier que
$F=\Q$. Un $\GL(d)$-torseur arithm\'etique $\hat E$ sur $\Spec\Z$
(pour le choix du sous-groupe compact maximal $U(d)$) n'est
autre qu'un $\Z$-module libre de rang~$d$ muni d'une norme euclidienne,
auquel on sait attacher (au moins) deux invariants :
$$ \theta(\hat E,t) = \sum_{e\in\hat E} \exp(-\pi t\norm{e}^ 2)
\quad\text{et}\quad
   \zeta(\hat E,s) = \sum_{e\in\hat E\setminus\{0\}}
                            \frac{1}{\norm{e}^ s}.       $$
(Ces s\'eries convergent respectivement pour $\Re(t)>0$ et $\Re(s)>d$.)
Comme il est bien connu,
la formule de Poisson standard implique l'\'equation fonctionnelle
$$ \theta(\hat E,t) = \frac{1}{t^{d/2}\vol(\hat E)} \theta(\hat E^ \vee,1/t) $$
o\`u $\vol(\hat E)=\exp(-\hdeg\hat E)$ est le covolume du r\'eseau
$\hat E$ dans $\hat E\otimes_\Z\R\simeq\R^ d$,
$\hat E^ \vee$ d\'esigne le r\'eseau dual (muni de la norme euclidienne duale)
et o\`u la d\'etermination de $t^ {d/2}$ est usuelle pour $t>0$.
Il est aussi bien connu comment utiliser cette \'equation pour en d\'eduire
que la fonction d\'efinie par 
$$ \Lambda(\hat E,s)
= \sqrt{\vol(\hat E)} \zeta(\hat E,s) \pi^ {-s/2}\Gamma(s/2) $$
poss\`ede un prolongement m\'eromorphe \`a $\C$, avec des p\^oles
simples en $s=0$ et $s=d$ de r\'esidus respectivement $-2\sqrt{\vol(\hat E)}$
et $2/\sqrt{\vol(\hat E)}$ et v\'erifie l'\'equation fonctionnelle
$$ \Lambda(\hat E,s)=\Lambda(\hat E^\vee,d-s). $$

Sur un corps de nombres quelconque, il faudrait tenir compte de la
diff\'erente, comme dans l'article r\'ecent
de van der Geer et Schoof~\cite{vdGeer-s98}. Selon ces m\^emes
auteurs, l'invariant $\theta(\hat E,1)$
mesure l'\emph{effectivit\'e} du fibr\'e vectoriel hermitien $\hat E$.
Ils interpr\`etent en particulier l'\'equation fonctionnelle de la 
fonction $\theta$ comme une formule de Riemann--Roch.

\paragraph{Exemples exotiques de fonctions $L$}
Soit maintenant $\hat {\mathcal E}\in\hatFib_d(\mathcal B)$.
On peut d\'efinir
des fonctions $L$ d'Arakelov
(pour une partie $U\subset\mathcal B(F)$ fix\'ee)
$$ \Theta(\hat {\mathcal E},s)
= L(\hat {\mathcal E},U, \theta(\cdot,1)\vol(\cdot)^s)
= \sum_{b\in U\subset \mathcal B(F)} \theta(\hat {\mathcal E}|_b,1)\vol(\hat {\mathcal E}|_b)^s $$
et
$$ Z(\hat {\mathcal E},s)= L(\hat {\mathcal E},U, \zeta(\cdot,ds)\vol(\cdot)^s)
= \sum_{b\in U\subset \mathcal B(F)} \zeta(\hat {\mathcal E}|_b,ds)\vol(\hat {\mathcal E}|_b)^s $$
et l'on a les \'egalit\'es, o\`u chacun des membres converge absolument
quand l'autre converge absolument,
$$ \Theta(\hat {\mathcal E},s)=\Theta(\hat {\mathcal E}^\vee,1-s)
\quad\text{et}\quad
    Z (\hat {\mathcal E},s) = Z(\hat {\mathcal E}^\vee,1-s). $$
Par exemple, pour $\mathcal B=\P^1_\Z$ et $\hat {\mathcal E}=\mathcal O_{\P}(1)$
avec la m\'etrique {\og max. des coordonn\'ees\fg}, on a
$$ \Theta(\hat {\mathcal E},s) = \sum_{N\geq 1} 2(1+2\phi(N)) \theta(N^2) N^{1-s}, $$
expression qui converge pour $\Re(s)>3$ et dans laquelle $\theta$
d\'esigne la fonction th\^eta de Riemann.

\section{Espaces fibr\'es}


\subsection{Constructions}
\label{constructions}

\begin{enonce}{Situation} \label{situation}
Soient $S$  un sch\'ema,
$G$ un $S$-sch\'ema en groupes lin\'eaire et plat,
dont on suppose pour simplifier les fibres g\'eom\'etriquement connexes
$f:\mathcal X\ra S$ un $S$-sch\'ema plat (quasi-compact et quasi-s\'epar\'e),
muni d'une action de $G/S$.
Soient aussi $g:\mathcal B\ra S$ un $S$-sch\'ema plat ainsi qu'un $G$-torseur
$\mathcal T\ra \mathcal B$ localement trivial pour la topologie
de Zariski.
\end{enonce}

\begin{enonce}{Construction}  \label{consY}
On d\'efinit un $S$-sch\'ema $\mathcal Y$, muni d'un morphisme
$\pi:\mathcal Y\ra \mathcal B$ localement isomorphe \`a $\mathcal X$
sur $\mathcal B$, par le changement de groupe structural
$G\ra \Aut_S(\mathcal X)$.
\end{enonce}
En effet, soit $(U_i)_{i\in I}$ un recouvrement ouvert
de $\mathcal B$ tel qu'il existe une trivialisation $\phi_i:G\times_S
U_i \xrightarrow{\sim} \mathcal T|_{U_i}$.
Si $i,j\in I$, soit $g_{ij}\in\Gamma(U_i\cap U_j,G)$
l'unique section telle que $\phi_i=g_{ij}\phi_j$ sur $U_i\cap U_j$.
En particulier, les $g_{ij}$ donnent un cocycle dont la classe
dans $\rmH^1(\mathcal B,G)$
repr\'esente la classe d'isomorphisme
du $G$-torseur $\mathcal T$.
Posons $\mathcal Y_i=\mathcal X\times_S U_i$ ; alors,
$g_{ij}$ agit sur $\mathcal X\times_S (U_i\cap U_j)$
et induit un isomorphisme
$$\phi_{ij}:\mathcal Y_j|_{U_i\cap U_j} \simeq \mathcal
Y_i|_{U_i\cap U_j}$$
que l'on utilise pour recoller les $\mathcal Y_i$.

On laisse v\'erifier que $\mathcal Y$ est un $\mathcal B$-sch\'ema bien d\'efini,
c'est-\`a-dire qu'il ne d\'epend pas \`a isomorphisme canonique pr\`es
du choix des trivialisations locales que l'on a fait.

\begin{lemm}
On a $\pi_*\mathcal O_{\mathcal Y}=g^*f_*\mathcal O_{\mathcal X}$.
\end{lemm}

\begin{rema}
Dans certains cas, $\mathcal Y$ h\'erite d'une action d'un sous-groupe
de $G$, notamment quand $G$ est commutatif.
\end{rema}

\begin{enonce}{Construction}\label{cons2}
Il r\'esulte de la construction pr\'ec\'edente une application
$$\vartheta:\mathrm Z^{d,G}(\mathcal X)\ra \mathrm Z^d(\mathcal Y)$$
des cycles $G$-invariants
de codimension $d$ sur $\mathcal X$ dans les cycles de codimension $d$
sur $\mathcal Y$.
\end{enonce}


\begin{defi}
Une \emph{$G$-lin\'earisation} d'un faisceau quasi-coh\'erent $\mathcal F$
sur $\mathcal X$ est une action de $G$ sur $\V(\mathcal F)$ qui rel\`eve
l'action de $G$ sur $\mathcal X$.

Un morphisme (resp.\ le produit tensoriel, le dual, la somme directe,
le faisceau des homomorphismes, des extensions, etc.) 
de faisceaux quasi-coh\'erents $G$-lin\'earis\'es est d\'efini
naturellement. On note $\Qcoh^G(\mathcal X)$
(resp.\ $\Fib_d^G(\mathcal X)$, resp.\ $\Pic^G(\mathcal X)$)
la cat\'egorie
des faisceaux quasi-coh\'erents
(resp. de fibr\'es vectoriels de rang $d$, resp.\ des classes d'isomorphisme
de fibr\'es inversibles) $G$-lin\'earis\'es sur $\mathcal X$.
\end{defi}

\begin{enonce}{Construction}\label{cons3}
On construit un foncteur
$$\vartheta:\Qcoh^G(\mathcal X)\ra\Qcoh(\mathcal Y)$$
qui est compatible avec les op\'erations standard sur les faisceaux
quasi-coh\'erents. 
\end{enonce}
Soit $\mathcal F$ un faisceau quasi-coh\'erent $G$-lin\'earis\'e sur $\mathcal
X$. Reprenons les notations de la construction~\ref{consY}
de $\mathcal Y$. Posons $\mathcal F_i$ le faisceau quasi-coh\'erent
sur $\mathcal Y_i=\mathcal X\times_S U_i$ image r\'eciproque de $\mathcal F$
par la premi\`ere projection. Gr\^ace \`a la $G$-lin\'earisation sur $\mathcal F$,
les $g_{ij}$ induisent des isomorphismes
$$\phi_{ij}^*\mathcal F_j|_{\mathcal X\times(U_i\cap U_j)}\simeq
\mathcal F_i|_{\mathcal X\times (U_i\cap U_j)}$$
qui fournissent par recollement un faisceau quasi-coh\'erent 
sur $\mathcal Y$.

On laisse v\'erifier que ce foncteur est bien d\'efini, c'est-\`a-dire,
est ind\'ependant des choix que l'on a fait.

Si $\mathcal F$ est un fibr\'e vectoriel $G$-lin\'earis\'e de rang $d$
sur $\mathcal X$, il est clair que le faisceau obtenu sur $\mathcal Y$
est aussi un fibr\'e vectoriel de rang $d$.

On laisse v\'erifier que cette application est compatible aux op\'erations
standard, et en particulier qu'elle descend en une application
sur les classes d'isomorphisme.

Un cas particulier des constructions pr\'ec\'edentes est obtenu lorsque
$\mathcal X=S$, auquel cas $\mathcal Y=\mathcal B$.
On notera $\eta_{\mathcal T}$ l'application qui en r\'esulte
des faisceaux quasi-coh\'erents sur $S$
avec action de $G/S$ vers les faisceaux quasi-coh\'erents
sur $\mathcal B$. 
Bien s\^ur, $\eta_{\mathcal T}:\Rep_d (G)\ra \Fib_d (\mathcal B)$
n'est autre que l'application
usuelle de changement de groupe structural (passage d'un $G$-torseur
\`a un $\GL(d)$-torseur). 

\begin{prop} \label{prop:canonique}
Le faisceau $\Omega^1_{\mathcal X/S}$ est muni
d'une lin\'earisation canonique de $G$.
Par la construction~\ref{cons3},
on obtient le faisceau $\Omega^1_{\mathcal X/\mathcal B}$.

Supposons en particulier que $\mathcal X$ et $\mathcal B$ sont lisses
sur $S$ ; le faisceau canonique sur $\mathcal X/S$ est alors
automatiquement $G$-lin\'earis\'e et on a un isomorphisme
$$\omega_{\mathcal Y/S} \simeq \vartheta(\omega_{\mathcal X/S})\otimes
\pi^*\omega_{\mathcal B/S}. $$
\end{prop}
\begin{proof}
Si $(U_i)$ est un recouvrement ouvert de $\mathcal B$ avec des
isomorphismes
$(\phi_i,\pi):\pi^{-1}(U_i)\simeq \mathcal X\times_S U_i$ comme dans la
construction~\ref{consY}, on a un isomorphisme naturel
$$
\Omega^1_{\mathcal Y/\mathcal B}|_{\pi^{-1}(U_i)} =
\Omega^1_{\mathcal \pi^{-1}(U_i)/U_i} 
\simeq \phi_i^*\Omega^1_{\mathcal X/S}
$$
qui se recollent pr\'ecis\'ement comme dans la construction~\ref{cons3}.

Dans le cas o\`u $\mathcal X/S$ et $\mathcal B/S$ sont lisses,
la suite exacte
$$ 0 \ra \Omega^1_{\mathcal Y/\mathcal B} \ra \Omega^1_{\mathcal Y/S}\ra
\Omega^1_{\mathcal B/S} \ra 0 $$
implique que
$$ \omega_{\mathcal Y/S}
 \simeq
  \det\Omega^1_{\mathcal Y/\mathcal B}\otimes \pi^*\omega_{\mathcal B/S} 
 \simeq
  \vartheta(\omega_{\mathcal X/\mathcal B})\otimes\pi^*\omega_{\mathcal B/S}.
\qquad\qed $$
\let\qed\relax
\end{proof}

\begin{lemm} \label{lemm:calcul1}
Si $\mathcal F\in \Qcoh^G(\mathcal X)$, $f_*\mathcal F$ est
muni d'une action naturelle de $G$ et
$g_*\vartheta(\mathcal F)$ est canoniquement isomorphe \`a
$\eta_{\mathcal T}(f_*\mathcal F)$.
\end{lemm}
\begin{proof} Laiss\'ee au lecteur.\end{proof}

\begin{prop}  \label{prop:effectif}
Soient $(\lambda,\alpha)\in\Pic^G(\mathcal X)\times\Pic(\mathcal B)$.
Le fibr\'e en droites $\vartheta(\lambda)\otimes\pi^*\alpha$ sur $\mathcal Y$
est effectif si et seulement si le fibr\'e vectoriel sur $\mathcal B$
$$ \eta_{\mathcal T}(f_*\lambda)\otimes\alpha $$
est effectif.
Cela implique que $\lambda$ est effectif.
\end{prop}
\begin{proof}
On a
$$ \pi_*(\vartheta(\lambda)\otimes\pi^*\alpha) 
= \pi_*(\vartheta(\lambda))\otimes\alpha 
= \eta_{\mathcal T}(f_*\lambda)\otimes\alpha $$
d'apr\`es le lemme~\ref{lemm:calcul1}.
\end{proof}

Notons $\iota$ le morphisme de groupes naturel
$ X^*(G)\ra\Pic^G(\mathcal X) $
qui associe \`a un caract\`ere $\chi$ le fibr\'e trivial 
muni de la lin\'earisation telle que $G$ agit par $\chi$
sur le second facteur de $\mathcal X\times_S \A^1_S$.
\begin{prop} \label{prop:iota-eta}
Pour tout caract\`ere $\chi$, il existe un isomorphisme canonique
de faisceaux inversibles
$$ \vartheta (\iota(\chi)) \simeq \pi^* \eta_{\mathcal T}(\chi). $$
\end{prop}
\begin{proof}
Soit  $(U_i)$ un recouvrement ouvert de $\mathcal B$
avec des isomorphismes
$(\phi_i,\pi):\pi^{-1}(U_i)\simeq\mathcal X\times_S U_i$ ;
notons $g_{ij}\in G(U_i\cap U_J)$ tel que
$\phi_i=g_{ij}\cdot \phi_j:\pi^{-1}(U_i\cap U_j)\ra \mathcal X$.
Alors, le fibr\'e en droites $\vartheta(\iota(\chi))$ est
obtenu en recollant $\A^ 1\times \mathcal X\times U_i$
et $\A^1\times\mathcal X\times U_j$ par le morphisme
$ (t,x,u)\mapsto (\chi(g_{ij})t,g_{ij}\cdot x, u)$.

D'autre part, $\eta_{\mathcal T}(\chi)$ est un fibr\'e
en droite sur $\mathcal B$ obtenu en recollant $\A^1\times U_i$
et $\A^ 1\times U_j$ par $(t,u)\mapsto (\chi(g_{ij})t,u)$.
\end{proof}

\subsection{Groupe de Picard}

Dans ce paragraphe, on suppose que $S$ est le spectre d'un corps $F$
de caract\'eristique~$0$.
On cherche \`a exprimer le groupe de Picard de $\mathcal Y$
en fonction de ceux de $\mathcal X$ et $\mathcal B$.
Pour cela, on se place sous les hypoth\`eses suivantes :
\paragraph{Hypoth\`eses sur $\mathcal X$}
\label{hyp:picard-X}
On suppose que 
\begin{enumerate}
\item $\mathcal X$ est propre, lisse, g\'eom\'etriquement int\`egre ;
\item $\rmH^1(\mathcal X,\mathcal O_{\mathcal X})=0$ ;
\item $\mathcal X(F)$ est non vide ;
\item tout fibr\'e en droites sur $\mathcal X$ est $G$-lin\'earisable,
et de m\^eme apr\`es toute extension alg\'ebrique de $F$ ;
\item $\Pic(\mathcal X_{\bar F})$ est sans torsion.
\end{enumerate}

\begin{rema}
Ces hypoth\`eses concernant $\mathcal X$ sont v\'erifi\'ees
lorsque $\mathcal X$ est une vari\'et\'e torique projective d\'eploy\'ee sur $F$,
ou bien un espace de drapeaux
g\'en\'eralis\'e pour un groupe alg\'ebrique d\'eploy\'e sur $F$.

Elles entra\^{\i}nent que les groupes de Picard et de N\'eron-S\'everi
de $\mathcal X_{\bar F}$ co\"{\i}ncident (voir la preuve du
lemme~\ref{lemm:picard} plus bas). En particulier,
$\Pic(\mathcal X_{\bar F})$ est sous ces hypoth\`eses un
$\Z$-module libre de rang fini.

D'autre part, il est prouv\'e dans~\cite{mumford-f-k94}, Cor.~1.6, p.~35,
que sous l'hypoth\`ese (i), tout fibr\'e en droites sur $\mathcal X$
admet une puissance $G$-lin\'earisable. (Rappelons que $G$ est
connexe.) Le lecteur qui d\'esirerait s'affranchir de cette hypoth\`ese
v\'erifiera que de nombreux r\'esultats de la suite de ce texte 
restent vrais, au moins apr\`es tensorisation par $\Q$.
\end{rema}


\begin{lemm}     \label{lemm:picard}
Si les hypoth\`eses~\ref{hyp:picard-X}
sont satisfaites, on a les deux assertions :
\begin{itemize}
\item $\rmH^0(\mathcal X,\mathcal O_{\mathcal X})=F$ ;
\item pour tout $F$-sch\'ema connexe $U$ poss\'edant un
point $F$-rationnel, l'homomorphisme
naturel
$$ \Pic(\mathcal X)\times \Pic(U) \ra \Pic(\mathcal X\times_F U)$$
est un isomorphisme.
\end{itemize}
\end{lemm}
\begin{proof}
La premi\`ere proposition d\'ecoule de la factorisation de Stein.
Pour la seconde, 
on a d'apr\`es~\cite[8.1/4]{bosch-l-r90} une suite exacte
$$ 0 \ra \Pic(U)\ra\Pic(\mathcal X\times_F U)\ra
    \Pic_{\mathcal X/F}(U) \ra 0 . $$
En particulier, $\Pic(\mathcal X)=\Pic_{\mathcal X/F}(F)$.
La nullit\'e de $\rmH^1(\mathcal X,\mathcal O_{\mathcal X})$ implique
que $\Pic_{\mathcal X/F}$ est de dimension~$0$,
donc que sa composante neutre $\Pic^0_{\mathcal X/F}=0$
puisque $F$ est de caract\'eristique nulle. Ainsi, $\Pic_{\mathcal X/F}$
est discret.
Alors, tout point rationnel $u\in U(F)$
d\'efinit un homomorphisme
$u^*:\Pic_{\mathcal X/F}(U)\ra\Pic_{\mathcal X/F}(F)$
qui par connexit\'e est l'inverse
de l'homomorphisme naturel
$\Pic_{\mathcal X/F}(F)\ra\Pic_{\mathcal X/F}(U)$.
\end{proof}

\begin{theo}\label{theo:picard}
Si $\iota$ d\'esigne le morphisme de groupes 
$ X^*(G)\ra\Pic^G(\mathcal X) $ introduit au paragraphe pr\'ec\'edent,
consid\'erons l'homomorphisme
$$
\Pic^G(\mathcal X)\oplus\Pic(\mathcal B) \ra \Pic(\mathcal Y), \quad
    (\lambda,\alpha) \mapsto \vartheta(\lambda)\otimes\pi^*\alpha.
$$
Si les hypoth\`eses~\ref{hyp:picard-X} sont satisfaites
et si $\mathcal B(F)$
est Zariski-dense dans~$\mathcal B$, 
alors la suite
$$ 0\ra X^*(G) 
\xrightarrow {(\iota,-\eta_{\mathcal T})}
       \Pic^G(\mathcal X)\oplus \Pic(\mathcal B)
\xrightarrow {\vartheta\otimes \pi^*} \Pic(\mathcal Y) \ra 0 $$
est exacte.
\end{theo}
\begin{proof}
Si $\iota(\chi)$ est trivial dans $\Pic^G(\mathcal X)$,
il r\'esulte de ce que $\rmH^0(\mathcal X,\mathcal O_{\mathcal X})=F$
que $\chi$ est n\'ecessairement le caract\`ere trivial.
En particulier, le premier homomorphisme est injectif.

La proposition~\ref{prop:iota-eta} implique que la composition
des deux premiers homomorphismes est nulle.

Si $\lambda$ est un fibr\'e en droites $G$-lin\'earis\'e sur $\mathcal X$
et $\alpha$ est un fibr\'e en droites sur $\mathcal B$, 
$\vartheta(\lambda)\otimes\pi^*\alpha$ est un fibr\'e en droites
sur $\mathcal Y$ dont la classe d'isomorphisme ne d\'epend que des
classes d'isomorphismes de $\lambda$ dans $\Pic^G(\mathcal X)$
et $\alpha$ dans $\Pic(\mathcal B)$.

Supposons qu'elle soit triviale. Soit $b$ un point
$F$-rationnel de $\mathcal B$.
En restreignant $\vartheta(\lambda)\otimes\pi^*\alpha$ \`a
$\pi^{-1}(b)$, la construction~\ref{cons2} de $\vartheta(\lambda)$
implique que $\lambda$ est trivial. 
La $G$-lin\'earisation de $\lambda$ est ainsi donn\'ee par un 
caract\`ere $\chi$ de $G$ et $\lambda=\iota(\chi)$.
D'apr\`es la proposition~\ref{prop:iota-eta},
on a $\vartheta(\lambda)=\pi^*\eta_{\mathcal T}(\chi)$.
Par suite, $\pi^*\alpha\simeq\pi^*\eta_{\mathcal T}(\chi)^{-1}$,
ce qui prouve l'exactitude au milieu.

Montrons alors que la derni\`ere fl\`eche est surjective.
Soit $\mathcal L$ un fibr\'e en droites sur $\mathcal Y$.
On peut recouvrir $\mathcal B$ par des ouverts connexes non vides $U_i$
assez petits de sorte que 
$$\pi^{-1}(U_i)\simeq \mathcal X\times_F U_i.$$
La restriction de $\mathcal L$ \`a $\pi^{-1}(U_i)$ fournit alors pour
tout~$i$ un \'el\'ement de
$$ \Pic(\mathcal X\times_F U_i)=\Pic(\mathcal X) \times \Pic(U_i) $$
puisque chaque $U_i$ a un point $F$-rationnel.
On en d\'eduit d'abord pour tout $i$ un \'el\'ement de $\Pic(\mathcal X)$ qui,
comme on le voit en les restreignant \`a $U_i\cap U_j$,
ne d\'epend pas de $i$. Notons le $\lambda$.
Finalement, il existe un faisceau inversible $\alpha_i\in\Pic(U_i)$
tel que la restriction de $\mathcal L$ \`a
$\pi^{-1}(U_i)\simeq \mathcal X\times_F U_i$ est
isomorphe \`a $p_1^*\lambda\otimes p_2^*\alpha_i$.
Quitte \`a raffiner le recouvrement $(U_i)$, on peut de plus supposer
que $\alpha_i\simeq \mathcal O_{U_i}$.

Choisissons une $G$-lin\'earisation sur $\lambda$.
On constate que la restriction de
$\mathcal L\otimes\vartheta(\lambda)^{-1}$
\`a $\pi^{-1}(U_i)$ est triviale. Si l'on choisit des trivialisations
on obtient en les comparant sur $\pi^{-1}(U_i\cap U_j)$ 
un \'el\'ement de 
$$\Gamma(\pi^{-1}(U_i\cap U_j), \mathcal O_{\mathcal Y}^\times)
     = \Gamma(U_i\cap U_j, \mathcal O_{\mathcal B}^\times) $$
car $\rmH^0(\mathcal X,\mathcal O_{\mathcal X})=F$.
Ces \'el\'ements d\'efinissent un $2$-cocycle de \v Cech sur $\mathcal B$
\`a valeurs dans le faisceau $\mathcal O_{\mathcal B}^\times$, 
d'o\`u un fibr\'e en droites $\alpha\in\Pic(\mathcal B)$ tel que
$$\mathcal L\otimes\vartheta(\lambda)^{-1} \simeq \pi^*\alpha. $$
Autrement dit, $\mathcal L$ appartient \`a l'image de l'homomorphisme
$\vartheta\otimes\pi^*$.

Le th\'eor\`eme est ainsi d\'emontr\'e.
\end{proof}

\begin{coro}
Supposons v\'erifi\'ees les hypoth\`eses~\ref{hyp:picard-X} 
et supposons que $\mathcal B(F)$ est Zariski-dense dans $\mathcal B$.
On dispose alors de suites exactes de $\Z[\Gal(\bar F/F)]$-modules :
\begin{gather}
0 \ra X^*(G_{\bar F}) \ra \Pic^G(\mathcal X_{\bar F})
  \ra \Pic(\mathcal X_{\bar F}) \ra 0  \\
0\ra X^*(G_{\bar F}) \ra
   \Pic^G(\mathcal X_{\bar F})\oplus\Pic(\mathcal B_ {\bar F}) \ra
   \Pic(\mathcal Y_{\bar F}) \ra 0 \\
 0 \ra \Pic(\mathcal B_{\bar F}) \xrightarrow{\pi^*}  
   \Pic (\mathcal Y_{\bar F}) \ra \Pic(\mathcal X_{\bar F})\ra 0. 
\end{gather}
\end{coro}
\begin{proof}
Il suffit d'appliquer le th\'eor\`eme~\ref{theo:picard}
sur $\bar F$, et de constater que la suite exacte obtenue
est $\Gal(\bar F/F)$-\'equivariante.
\end{proof}

\begin{theo}  \label{theo:effectif}
Supposons v\'erifi\'ees les hypoth\`eses~\ref{hyp:picard-X},
que $\mathcal B(F)$ est Zariski-dense dans $\mathcal B$,
et supposons de plus que $G$ est un groupe alg\'ebrique $F$-r\'esoluble\footnote
{Cela signifie que $G$ est extension it\'er\'ee de $\gm$ et $\ga$,
autrement dit, que $G$ est r\'esoluble et d\'eploy\'e sur $F$.},
un fibr\'e en droites sur $\mathcal Y$ est alors effectif
si et seulement s'il s'\'ecrit comme l'image 
d'un couple $(\lambda,\alpha)\in\Pic^G(\mathcal X)\times\Pic(\mathcal B)$
o\`u $\lambda$ et $\alpha$ sont effectifs.
\end{theo}
\begin{proof}
Soient $\lambda\in\Pic^G(\mathcal X)$ et $\alpha\in\Pic(\mathcal B)$
effectifs. On veut montrer que $\vartheta(\lambda)\otimes\pi^*\alpha$
est effectif. Il suffit de prouver que $\vartheta(\lambda)$ est
effectif, et pour cela, il suffit de prouver qu'il existe un
diviseur de Cartier $G$-invariant $D$ sur $\mathcal X$
tel que l'on ait un isomorphisme de fibr\'es
en droites $G$-lin\'earis\'es, $\lambda\simeq\mathcal O(D)$.
Autrement dit, il faut montrer que la repr\'esentation de $G$
sur $f_*\lambda$ admet une $F$-droite stable, ce qu'implique
le th\'eor\`eme de point fixe de Borel puisque $G$ est $F$-r\'esoluble.

Soit maintenant $\mathcal L$ un fibr\'e en droites effectif sur $\mathcal Y$.
Comme $G$ est connexe et $\Pic(G)=0$,
la d\'emonstration de la proposition 1.5, p.~34, de~\cite{mumford-f-k94}
implique que tout fibr\'e inversible sur $\mathcal X$ est $G$-lin\'earisable.
Le th\'eor\`eme~\ref{theo:picard} implique donc qu'il existe
$\lambda\in\Pic^G(\mathcal X)$ et $\alpha\in\Pic(\mathcal B)$
tels que $\mathcal L=\vartheta(\lambda)\otimes\pi^*\alpha$.
D'apr\`es la proposition~\ref{prop:effectif},
$\eta_{\mathcal T}(f_*\lambda)\otimes\alpha$ est effectif.
Comme $G$ est $F$-r\'esoluble, toute repr\'esentation lin\'eaire de $G$
est extension successive de repr\'esentations de dimension~$1$.
Cela implique que $\eta_{\mathcal T}(f_*\lambda)$ est extension successive
de fibr\'es en droites ; notons les $\lambda_i$. Alors, $\eta_{\mathcal
T}(f_*\lambda)\otimes\alpha$ est extension des $\lambda_i\otimes\alpha$,
et l'effectivit\'e de $\mathcal L$ implique que l'un au moins
des $\lambda_i\otimes\alpha$ est effectif.

Or, $\lambda_i$ est associ\'e \`a un caract\`ere $\chi_i$ de $G$ ;
si on remplace $\lambda$ par le fibr\'e en droite
$G$-lin\'earis\'e $\lambda\otimes \iota (\chi_i)^{-1}$
o\`u l'action a \'et\'e divis\'ee par $\chi_i$, on repr\'esente ainsi
$\mathcal L$ sous la forme
$$ \mathcal L \simeq \vartheta(\lambda\otimes \iota(\chi_i)^{-1}) \otimes
(\lambda_i\otimes\alpha) , $$
ce qui conclut la d\'emonstration, $\lambda\otimes\iota(\chi_i)^{-1}$
\'etant isomorphe \`a $\lambda$ comme fibr\'e en droites, donc effectif.
\end{proof}


\subsection{M\'etriques hermitiennes}

Dans ce paragraphe, nous \'etendons la construction~\ref{cons3}
en supposant que $S$ est le spectre d'un corps de nombres
et en faisant intervenir des m\'etriques hermitiennes.

\begin{defi}
Soit $G$ un groupe de Lie connexe sur $\C$ ; fixons un sous-groupe
compact maximal $K$ de $G$.
Soit $X$ une vari\'et\'e analytique complexe munie d'une action de $G$.

Si $\mathcal E$ est un fibr\'e vectoriel complexe $G$-lin\'earis\'e
sur $X$, on dit qu'une m\'etrique hermitienne est \emph{$K$-invariante}
si l'action de $K$ sur $\V(\mathcal E)\times X$ est isom\'etrique.
\end{defi}

On remarquera que les constructions usuelles (tensorielles) de
fibr\'es hermitiens pr\'eservent la $K$-invariance
des m\'etriques hermitiennes.

\begin{rema}
Avec les notations de la d\'efinition pr\'ec\'edente, tout fibr\'e vectoriel
sur $X$ admet une m\'etrique hermitienne $K$-invariante :
si $\norm{\cdot}_0$ est une m\'etrique hermitienne sur $\mathcal E$,
on peut en effet choisir une mesure de Haar sur $K$ et poser
pour toute section $\mathsf s$,
$$\norm{\mathsf s}^2(x) = \int_{K} \norm{k\cdot \mathsf s}(x) ^2\, dk .$$
\end{rema}

Rappelons l'\'enonc\'e de la situation~\ref{situ:GK} :
\begin{enonce*}{Situation} 
Supposons que $S$ est le spectre de l'anneau des entiers d'un corps de nombres
$F$ et que $G$ est un $S$-sch\'ema en groupes lin\'eaire connexe.
Fixons pour tout plongement complexe de $F$ $\sigma\in S(\C)$ un sous-groupe
compact maximal $K_\sigma$ de $G(\C)$ et notons
$\mathbf K_\infty$ la collection $(K_\sigma)_{\sigma}$.
\end{enonce*}

\begin{defi}
Supposons que $G$ agit sur un $S$-sch\'ema plat $\mathcal X$.
On appelle \emph{fibr\'e vectoriel hermitien $(G,\mathbf K_\infty)$-lin\'earis\'e} un
fibr\'e vectoriel $\mathcal E$ sur $\mathcal X$ muni d'une $G$-lin\'earisation
et, pour tout $\sigma\in S(\C)$,
d'une m\'etrique hermitienne sur le fibr\'e vectoriel
$\mathcal E\otimes_\sigma\C$ sur $\mathcal X(\C)$ qui est
$K_\sigma$-invariante.

On note $\hatFib^{G,\mathbf K_\infty}_d(\mathcal X)$ la cat\'egorie
des fibr\'es vectoriels hermitiens $(G,\mathbf K_\infty)$-lin\'earis\'es
de rang~$d$ sur $\mathcal X$.
Si $d=1$, on notera $\hatPic^{G,\mathbf K_\infty}(\mathcal X)$ le groupe
des classes d'isomorphisme de fibr\'es vectoriels hermitiens de rang~$1$
$(G,\mathbf K_\infty)$-lin\'earis\'es sur $\mathcal X$.
\end{defi}

\begin{enonce}{Situation}\label{hatsitu}
Pla\c cons-nous dans la situation~\ref{situ:GK}.
Soit $f:\mathcal X\ra S$ un $S$-sch\'ema plat,
muni d'une action de $G/S$.
Soient aussi $g:\mathcal B\ra S$ un $S$-sch\'ema plat ainsi qu'un $(G,\mathbf
K)$-torseur arithm\'etique $\hat{\mathcal T}$ sur $\mathcal B$
(voir la d\'efinition~\ref{defi:torsarith}).
\end{enonce}

%

\begin{enonce}{Construction}
\label{cons:hat}
Le foncteur $\vartheta:\Fib_d^G(\mathcal X)\ra \Fib_d(\mathcal Y)$ s'\'etend
en un foncteur
$$ \vartheta: \hatFib_d^{G,\mathbf K_\infty}(\mathcal X)
\ra \hatFib_d(\mathcal Y) $$
qui est compatible avec les op\'erations tensorielles
standard sur les fibr\'es vectoriels
hermitiens $(G,\mathbf K_\infty)$-lin\'earis\'es
(resp.\  les fibr\'es vectoriels hermitiens).
\end{enonce}
Soit $\mathcal F$ un fibr\'e vectoriel hermitien
$(G, \mathbf K_\infty)$-lin\'earis\'e
sur $\mathcal X$.
Soit $\sigma\in S(\C)$. 
De mani\`ere analogue \`a ce qu'on a fait
dans la construction~\ref{cons2}, choisissons 
un recouvrement ouvert $(U_i)$ de $\mathcal B_\sigma(\C)$
pour la topologie
complexe de sorte que la restriction du torseur $\mathcal T$ \`a $U_i$
est triviale et qu'il existe des trivialisations dont les fonctions
de transistions associ\'es $g_{ij}\in \Gamma(U_i\cap U_j,G)$ soient \`a
valeurs dans $K_\sigma$.
Le choix de telles trivialisations induit des isomorphismes
$$\pi^{-1}(U_i)\simeq \mathcal X(\C)\times U_i, \qquad
   \vartheta(\mathcal F)|_{\pi^{-1}(U_i)} \simeq p_1^*\mathcal F. $$
Pour tout~$i$,
on a ainsi une m\'etrique hermitienne naturelle sur $\vartheta(\mathcal
F)|_{\pi^{-1}(U_i)}$ par image r\'eciproque de la m\'etrique hermitienne
sur $\mathcal F$. Comme $g_{ij}\in K_\sigma$ et comme la m\'etrique
hermitienne sur $\mathcal F$ est $K_\sigma$-invariante, les
m\'etriques hermitiennes sur $\vartheta(\mathcal F)|_{U_i\cap U_j}$ induites
par $U_i$ et par $U_j$ co\"\i ncident, d'o\`u une m\'etrique hermitienne
bien d\'efinie sur $\vartheta(\mathcal F)$.

Enfin, la proposition~\ref{prop:iota-eta} admet une g\'en\'eralisation
avec m\'etriques hermitiennes :
\begin{prop} \label{prop:hat-iota-eta}
Pour tout caract\`ere $\chi\in X^*(G)$, l'isomorphisme
canonique de la proposition~\ref{prop:iota-eta} est
une isom\'etrie.
\end{prop}
\begin{proof}
Si l'on reproduit la d\'emonstration
de la proposition~\ref{prop:iota-eta} pour un recouvrement ouvert
pour la topologie complexe (les $g_{ij}$ \'etant donc dans
le sous-groupe compact maximal), chacun des fibr\'es est d\'efini
par recollement de la m\^eme mani\`ere, et les m\'etriques
sur ces fibr\'es sont d\'efinies de sorte que cette identification
soit une isom\'etrie. Il en r\'esulte que l'isomorphisme
de cette proposition,
qui consistait en l'application \'evidente sur les ouverts
$\mathcal X\times U_i$ est une isom\'etrie.
\end{proof}

\subsection{Torsion des m\'etriques ad\'eliques} \label{subsec:torsion}

Pla\c cons nous alors dans la situation~\ref{hatsitu}, toujours
avec $S=\Spec\mathfrak o_F$.
Soit $\mathcal L$ un fibr\'e en droites hermitien
$(G,\mathbf K_\infty)$-lin\'earis\'e sur $\mathcal X$. La restriction
de $\mathcal L$ \`a $\mathcal X_ F$ est ainsi
munie d'une m\'etrique ad\'elique naturelle.

\begin{enonce}{Proposition-D\'efinition}
Soit $\mathbf g=(g_v)_v\in G(\A_F)$. On d\'efinit une m\'etrique ad\'elique sur
$\mathcal L$, appel\'ee {\em m\'etrique ad\'elique tordue par $\mathbf g$}
en posant pour toute place $v$ de $F$, tout point $x\in \mathcal X(F_v)$
et toute section $\mathsf s\in\mathcal L_x$,
$$ \norm{\mathsf s}_v' (x) = \norm{g_v\cdot \mathsf s}_v (g_v\cdot x). $$
\end{enonce}
\begin{proof}
Il est clair que pour toute place $v$, on a d\'efini une m\'etrique $v$-adique.
L'ensemble des places non-archim\'ediennes $v$ telles que
$g_v\in G(\mathfrak o_v)$ est par d\'efinition de compl\'ementaire fini. 
Pour ces places, $\norm{\mathsf s}_v' (x)=\norm{\mathsf s}_v(x)$ car $g_v$
\'etant un automorphisme de $\mathcal L$ sur $\Spec\mathfrak o_v$,
la section $g_v\cdot s$ est enti\`ere en $g_v\cdot x$ si et seulement
si la section $s$ est enti\`ere en $x$.
Ainsi, hors d'un nombre fini de places, la nouvelle collection
de m\'etriques $v$-adiques est d\'efinie par un mod\`ele entier.
Elle d\'efinit donc une m\'etrique ad\'elique.
\end{proof}
Remarquons que $G(\A_F)$ n'agit en fait qu'\`a travers
$ G(\A_F)/ \mathbf K_G$.

\begin{exem}
Soit $E$ un $F$-espace vectoriel de dimension finie et notons
$\P$ l'espace projectif des droites de $E$. Faisons agir
$\GL(E)$ de mani\`ere naturelle sur $\P$.
Le faisceau $\mathcal O_{\P}(1)$
poss\`ede une $\GL(E)$-lin\'earisation naturelle d\`es qu'on
a remarqu\'e qu'une section
de $\mathcal O_{\P}(-1)$ en un point $\mathbf x\in \P$ correspond \`a un point
de la droite $D_{\mathbf x}$ d\'efinie par $\mathbf x$.
De mani\`ere explicite, l'espace vectoriel des sections globales
de $\mathcal O(1)$ sur $\P$ s'identifie au dual $E^*$ de $E$
sur lequel la $\GL(E)$-lin\'earisation sur $\mathcal O(1)$
induit la repr\'esentation contragr\'ediente $\phi\mapsto \phi\circ g^{-1}$.

Supposons que $E$ est muni d'une m\'etrique ad\'elique. On a alors une
m\'etrique ad\'elique sur $\mathcal O_{\P}(1)$ par la formule
$$ \norm{\phi} (\mathbf x) = \frac{\abs{\phi(e)}}{\norm{e}_v},
  \quad \phi\in E^*, \quad e \in D_{\mathbf x}\setminus\{0\}. $$
Il r\'esulte de la formule du produit
que la hauteur exponentielle d'un point $\mathbf x\in \P(F)$
est donn\'ee par la formule
$$ H(\mathbf x) = \prod_{v} \norm{e}_v, \quad e \in D_{\mathbf x}\setminus\{0\}
.$$

Soit alors $(g_v)_v\in \GL(E)(\A_F)$. La m\'etrique $v$-adique
tordue par $g_v$ sur $\mathcal O_{\P}(1)$ est ainsi donn\'ee par
$$ \norm{\phi}' (\mathbf x) = \frac{\abs{\phi(e)}}{\norm{g_v\cdot e}_v},
  \quad \phi\in E^*, \quad e \in D_{\mathbf x}\setminus\{0\}. $$
Autrement dit, la hauteur exponentielle tordue de $\mathbf x\in \P(F)$
est d\'efinie par l'expression
$$ H'(\mathbf x) = \prod_{v} \norm{g_v\cdot e}_v,
\quad e \in D_{\mathbf x}\setminus\{0\} .$$
Cette formule \'etait donn\'ee comme d\'efinition de la hauteur tordue
par Roy et Thunder dans~\cite{roy-t96}.
\end{exem}

Dans certains cas, on peut comparer la m\'etrique ad\'elique initiale
sur $\mathcal L$ et la m\'etrique ad\'elique tordue.

\begin{prop} \label{prop:comp-tordue}
Supposons que $\mathsf s$ est une section globale de $\mathcal L$
sur $\mathcal X_ F$
dont le diviseur est $G$-invariant.
Il existe alors un unique caract\`ere $\chi\in X^*(G)$ $F$-rationnel
(le \emph{poids} de $\mathsf s$) tel que pour tout $g\in G$,
$g\cdot \mathsf s = \chi(g) \mathsf s$.

Soit $\mathbf g\in G(\A_F)$, et consid\'erons $\overline{\mathcal L'}$
la m\'etrique ad\'elique tordue par $\mathbf g$.
Si $x\in X(F)$ n'appartient pas au diviseur de $\mathsf s$,
on a la formule
$$ H(\overline{\mathcal L'};x) = \prod_v \abs{\chi(g_v)}_v^{-1} \,
    H (\overline{\mathcal L}, \mathsf s,  \mathbf g\cdot x). $$
\end{prop}
\begin{proof}
Comme le diviseur de $\mathsf s$ est $G$-invariant, 
il existe pour tout $g\in G$ un \'el\'ement $\chi(g)\neq 0$
tel que $g\cdot \mathsf s=\chi(g)\mathsf s$.
Il est alors clair que $g\mapsto \chi(g)$ d\'efinit un caract\`ere $F$-rationnel
(alg\'ebrique) de $g$.

D'autre part, on a  pour toute place $v$ de $F$,
$$ \norm{\mathsf s}_v' (x) = \norm{g_v\cdot \mathsf s}_v (g_v x)
   = \norm{\chi(g_v)\mathsf s}_v (g_v x)
   = \abs{\chi(g_v)}_v \norm{\mathsf s}_v (g_v x). $$
La proposition en d\'ecoule en prenant le produit.
\end{proof}

\begin{rema}
Bien s\^ur, dans l'\'enonc\'e pr\'ec\'edent, il suffit de supposer
que la section $\mathsf s$ est propre pour les \'el\'ements $g_v$.
En particulier, si $G'$ est un sous-groupe de $G$ tel que
$\div(\mathsf s)$ est invariant par $G'$, on aura une formule
du m\^eme type pour les m\'etriques ad\'eliques tordue par un \'el\'ement
de $G(\A_F)$.
\end{rema}

\begin{rema}[Choix des sections]
La formule pr\'ec\'edente permet de comparer la restriction
\`a $G(\A_F)\mathcal X(F)$ des hauteurs sur les points ad\'eliques
associ\'ees \`a deux sections $\mathsf s_1$
et $\mathsf s_2$ de poids respectivement $\chi_1$ et $\chi_2$.
En effet, si $\mathbf x=\mathbf g\cdot x\in G(\A_F)\mathcal X(F)$,
on a, $\overline{\mathcal L'}$ d\'esignant la m\'etrique ad\'elique
tordue par $\mathbf g$,
$$ H(\overline{\mathcal L},\mathsf s_1; \mathbf x)
= \prod_v \abs{\chi_1(g_v)}_v \, H(\overline{\mathcal L'};x)
= \prod_v \abs{\chi_1\chi_2^{-1}(g_v)}_v \, H(\overline{\mathcal L},\mathsf
s_2; \mathbf x).
$$
Appliqu\'ee \`a des sections de m\^eme poids $\chi$, cela permet d'\'etendre
les fonctions $H(\overline{\mathcal L},\mathsf s;\cdot)$
au compl\'ementaire dans $G(\A_F)\mathcal X(F)$ de l'intersection
des diviseurs des sections de poids~$\chi$.
\end{rema}

\begin{rema} \label{rema:torsion-torique}
Lorsque $\mathcal X$ est une vari\'et\'e torique, compactification
\'equivariante lisse d'un tore $G$, tout fibr\'e en droites  effectif
$\mathcal L$ qui est $G$-lin\'earis\'e 
poss\`ede une unique droite $F$-rationnelle
de sections pour lesquelles $G$ agit par le caract\`ere trivial.
On peut utiliser cette section pour d\'efinir une hauteur
sur les points ad\'eliques du compl\'ementaire de son diviseur,
donc en particulier sur $G(\A_F)$.
\end{rema}

Expliquons maintenant comment la torsion des m\'etriques ad\'eliques
intervient dans nos constructions.
Nous allons pr\'eciser un peu la situation~\ref{hatsitu} en faisant
d\'esormais l'hypoth\`ese suivante :
\begin{enonce}{Situation}\label{hautsitu}
Nous faisons les hypoth\`eses contenues dans la situation~\ref{hatsitu}.
En particulier, $S$ est le spectre de l'anneau des entiers de corps
de nombres~$F$.
De plus, supposons que $\mathcal X$ et $\mathcal B$ sont propres
sur $S$.
\end{enonce}

Soit $b$ un point $F$-rationnel de $\mathcal B$. Comme $\mathcal B$
est propre sur $S$, il en r\'esulte une unique section $\eps_b:S\ra\mathcal
B$ qui prolonge $b$.
Toute trivialisation du $G_ F$-torseur
$G_ F\simeq \mathcal T|_b$ sur $\Spec F$
(il en existe car c'est un torseur pour la topologie de Zariski)
induit un isomorphisme $\mathcal X _ F\simeq \mathcal Y|_b$.
Fixons un tel isomorphisme $\phi$. Si $\lambda\in\Pic^G(\mathcal X)$,
$\phi^*\vartheta(\lambda)$ est un fibr\'e en droite sur $\mathcal X_ F$
canoniquement isomorphe \`a $\lambda$.
En revanche, les m\'etriques (ad\'eliques) sont en g\'en\'eral distinctes.

Soit $v$ une place finie de $F$, notons $\mathfrak o_v$
le compl\'et\'e de l'anneau local de $\mathfrak o_F$ en $v$.
Soit $\eps_v:\Spec\mathfrak o_v\ra \mathcal B$ la restriction
de $\eps_b$ \`a $\Spec\mathfrak o_v$.
Alors, $\eps_v^*\mathcal T$ est un $G\otimes \mathfrak o_v$-torseur
sur $\Spec\mathfrak o_v$, et est donc trivialisable.
Ainsi, $\eps_v^*\mathcal Y$ est isomorphe \`a $\mathcal X\otimes\mathfrak o_v$.
Fixons un isomorphisme $\phi_v$ induit par une trivialisation du
torseur. Il existe par d\'efinition $g_v\in G(F_v)$ tel que
$$\phi=\phi_v\circ [g_v],
    \quad \mathcal X\otimes F_v\ra \mathcal Y|_b\otimes F_v,
$$
$[g_v]$ d\'esignant l'automorphisme de $\mathcal X\otimes F_v$
d\'efini par $g_v$.
La d\'efinition de la m\'etrique $v$-adique associ\'ee \`a un mod\`ele
montre que $\phi_v$ est une isom\'etrie. Ainsi, en tant
que fibr\'e inversible m\'etris\'e sur $\mathcal X\otimes F_v$,
$ \phi^*(\vartheta(\lambda)) $ est isomorphe (isom\'etrique) \`a  $[g_v]^*\lambda$.

Soit maintenant $v$ une place \`a l'infini. Comme on s'\'etait fix\'e une
trivialisation du $G(\C)/K_v$-fibr\'e sur $\mathcal B(\C)$,
on dispose d'un isomorphisme $\phi_v$ bien d\'efini modulo $K_v$
qui par d\'efinition ne modifie pas les m\'etriques. La comparaison
entre $\phi$ et $\phi_v$ se fait comme pr\'ec\'edemment par
un \'el\'ement $g_v\in G(\C)$.

Il en r\'esulte le th\'eor\`eme :
\begin{theo} \label{theo:comp-torsarith}
Soit $\mathbf g=(g_v)_v\in G(\A_F)$ l'\'el\'ement du groupe ad\'elique
que nous venons d'introduire. Il repr\'esente
la classe de la restriction \`a $b$ du $(G,\mathbf K_\infty)$-torseur 
arithm\'etique $\hat{\mathcal T}$
dans l'isomorphisme de la proposition~\ref{prop:H1-adelique}.
De plus, la m\'etrique ad\'elique image
r\'eciproque sur $\phi^*\vartheta(\lambda)$ s'identifie \`a la m\'etrique
ad\'elique tordue par $\mathbf g$ sur $\lambda$.
\end{theo}

\subsection{Nombres de Tamagawa}

Commen\c cons par rappeler la d\'efinition, due \`a Peyre (cf.~\cite{peyre95}
et~\cite{peyre98})
des nombres de Tamagawa associ\'es \`a une m\'etrique ad\'elique
sur le faisceau anticanonique.

\paragraph{Hypoth\`eses} \label{tama:hypotheses}
Soit $X$ une vari\'et\'e propre, lisse et g\'eom\'etriquement int\`egre sur $F$
telle que $\rmH^1(X,\mathcal O_X)=\rmH^2(X,\mathcal O_X)=0$
et que $X(F)$ soit Zariski-dense dans $X$.
Sous ces conditions, $\Pic(X_{\bar F})_\Q$ est un $\Q$-espace
vectoriel de dimension finie.

\paragraph{D\'efinition}
Munissons le fibr\'e canonique $\omega_X$ d'une m\'etrique ad\'elique.
Pour toute place $v$ de $F$, une construction classique de Weil
fournit une mesure $\mu_{X,v}$ sur $X(F_v)$ \`a partir de la m\'etrique
$v$-adique sur $\omega_X$. Notons $L_v(s,\Pic(X_{\bar F}))$ 
le facteur local en $v$ de la fonction $L$ d'Artin de la repr\'esentation
de $\Gal(\bar F/F)$ sur $\Pic(X_{\bar F})_\Q$.
Le th\'eor\`eme de Weil sur la mesure de $X(F_v)$ pour $\mu_{X,v}$ 
et le th\'eor\`eme de Deligne sur les conjectures de Weil concernant
le nombre de points rationnels des vari\'et\'es sur les corps finis
ont la cons\'equence suivante : il existe un ensemble fini $\Sigma$ de places
de $F$, contenant les places archim\'ediennes, tel que
$$
\prod_{v\in\Sigma} \mu_{X,v}
\times
\prod_{v\not\in \Sigma}  \left( L^{-1}_v (1,\Pic(X_{\bar F})) \mu_{X,v} \right)
$$
d\'efinisse une mesure $\mu_{X,\Sigma}$ sur $X(\A_F)$ pour laquelle
$X(\A_F)$ a un volume fini.

Soit $L_\Sigma(s,\Pic(X_{\bar F}))
= \prod_{v\not\in\Sigma} L_v(s,\Pic(X_{\bar F}))$
la fonction $L$ partielle de $\Pic(X_{\bar F})$. Le produit eul\'erien
converge en effet pour $\Re(s)>1$ et $L_\Sigma$ a un p\^ole
en $s=1$ d'ordre
la dimension $t$ des invariants sous $\Gal(\bar F/F)$ de $\Pic(X_{\bar F})_\Q$.
Notons 
$$L^*_\Sigma(1,\Pic(X_{\bar F}))= \lim_{s\ra 1} (s-1)^r L_\Sigma(s,\Pic(X_{\bar F})). $$
On d\'efinit alors le nombre de Tamagawa de $X$ (associ\'e \`a la m\'etrique
ad\'elique choisie sur $\omega_X$) par
$$ \tau(X) = L^*_\Sigma(s,\Pic(X_{\bar F})) \,
        \int_{\overline{X(F)}} \mu_{X,\Sigma}. $$
Il est facile de v\'erifier qu'il ne d\'epend pas de l'ensemble fini
de places $\Sigma$ choisi.

Nous aurons \`a utiliser le lemme suivant.
\begin{lemm} \label{lemm:tamaouvert}
Supposons r\'ealis\'ees les hypoth\`eses~\ref{tama:hypotheses}.
Soit $U$ un ouvert non vide de $X$. 
Notons $\overline{U(F)}$ l'adh\'erence de $U(F)$ dans $\prod_v U(F_v)$
pour la topologie produit (qui est la topologie induite
sur $\prod_v U(F_v)$ par la topologie ad\'elique de $X(\A_F)$).
Alors, on a l'\'egalit\'e
$$ \int_{\overline{U(F)}} \mu_{X,\Sigma} =
\int_{\overline{X(F)}}\mu_{X,\Sigma}. $$
\end{lemm}
\begin{proof}
Tout point $x=(x_v)\in \prod_v U(F_v)$ poss\`ede par d\'efinition un
voisinage (pour la topologie induite) contenu dans $\prod_v U(F_v)$.
Par suite, si une suite $(x^{(n)})$ de points de $X(F)$ converge 
vers $x$, \`a partir d'un certain rang, $x^{(n)}$ appartient \`a $U(F_v)$
pour toute place $v$, et donc $x^{(n)}\in U(F)$. Cela montre
que $\overline{U(F)}=\overline{X(F)}\cap \prod_v U(F_v)$.
Ainsi, le compl\'ementaire de $\overline{U(F)}$ dans $\overline{X(F)}$
est contenu dans $X(\A_F)\setminus \prod_v U(F_v)$, donc dans la r\'eunion
$$ \bigcup_v (X\setminus U)(F_v) \prod_{w\neq v} X(F_w) . $$
La d\'efinition de la mesure $\mu_{X,v}$ implique que
$(X\setminus U)(F_v)$ est de mesure nulle pour $\mu_{X,v}$. 
On voit donc que $\overline{X(F)}\setminus \overline{U(F)}$
est r\'eunion d\'enombrable d'ensembles de mesure
nulle pour la mesure de Tamagawa sur $X(\A_F)$, donc est de mesure nulle.
\end{proof}

\medskip
On se place maintenant dans la situation~\ref{hatsitu},
$S$ \'etant le spectre $\Spec\mathfrak o_F$ de l'anneau des entiers 
d'un corps de nombres $F$.

\begin{lemm}
Si $\mathcal X_ F$ et $\mathcal B_ F$
satisfont les hypoth\`eses~\ref{tama:hypotheses} n\'ecessaires
pour la d\'efinition des nombres de Tamagawa,
$\mathcal Y_ F$ les satisfait aussi.
\end{lemm}
\begin{proof}
Que $\mathcal Y_F$ soit lisse, propre et g\'eom\'etrique int\`egre est clair.
D'autre part, les points rationnels de $\mathcal Y_F$ sont
denses dans chaque fibre au-dessus d'un point rationnel de $\mathcal B_F$,
lesquels sont suppos\'es denses dans $\mathcal B_F$. Comme
$\mathcal Y_F\ra\mathcal B_F$ est propre, un argument \'el\'ementaire
de platitude puis de dimension implique que
les points rationnels de $\mathcal Y_F$ sont Zariski-denses.

D'autre part, les hypoth\`eses sur $\mathcal X_F$ impliquent
que $R^0\pi_*\mathcal O_{\mathcal Y_F}=\mathcal O_{\mathcal B_F}$
et que 
$$
R^1\pi_* \mathcal O_{\mathcal Y_F}=R^2\pi_* \mathcal O_{\mathcal Y_F}=0.
$$
La suite spectrale des foncteurs
compos\'es implique que $\rmH^{q}(\mathcal O_{\mathcal B_F})$
est un quotient de $\bigoplus_{i+j=q}
\rmH^i(\mathcal B_F,R^j\pi_*\mathcal O_{\mathcal Y_F})$.
Si $j=1$ ou si $j=2$,
on a $\rmH^i(R^j\pi_*)=0$ puisque $R^j\pi_*=0$.
Si $j=0$ et $i\in\{1,2\}$, $\rmH^i(R^0\pi_*)
=\rmH^i(\mathcal O_{\mathcal B_F})=0$
en vertu des hypoth\`eses faites sur $\mathcal B_F$.
\end{proof}

Supposons donc que  $\mathcal X_ F$ et $\mathcal B_ F$
satisfont ces hypoth\`eses~\ref{tama:hypotheses}.
Le faisceau canonique sur $\mathcal Y$ admet d'apr\`es la
proposition~\ref{prop:canonique} une d\'ecomposition 
$$ \omega_{\mathcal Y}=
\vartheta(\omega_{\mathcal X}/S)\otimes\pi^*\omega_{\mathcal B/S}. $$
Choisissons une structure de fibr\'e en droite hermitien
$(G,\mathbf K_\infty)$ lin\'earis\'e sur $\omega_{\mathcal X/S}$
compatible \`a la lin\'earisation canonique sur $\omega_{\mathcal X/S}$
(autrement dit, pour toute place archim\'edienne $\sigma$,
une m\'etrique hermitienne $K_\sigma$-invariante sur $\mathcal
X\times_\sigma\C$).
Choisissons aussi une m\'etrique hermitienne sur $\omega_{\mathcal B/S}$.
Il en r\'esulte une m\'etrique hermitienne canonique sur $\omega_{\mathcal
Y/S}$ par la construction~\ref{cons:hat}.
Le fait de disposer d'un mod\`ele sur $\mathfrak o_F$ induit
de plus des m\'etriques $v$-adiques au places finies,
d'o\`u des m\'etriques ad\'eliques sur $\omega_{\mathcal X_ F}$,
sur $\omega_{\mathcal B_ F}$ et sur $\omega_{\mathcal Y_ F}$.

\begin{theo}
Muni de ces m\'etriques ad\'eliques, on a l'\'egalit\'e
$$ \tau(\mathcal Y_ F)
    = \tau(\mathcal X_ F)\tau(\mathcal B_ F). $$
\end{theo}

\begin{proof}
Soit $U$ un ouvert de Zariski non vide de $\mathcal B_F$
tel que $\mathcal T|_U\simeq G\times_S U$. Notons
$V=\pi^{-1}(U)\subset\mathcal Y_F$, de sorte que $V$
est un ouvert non vide de $\mathcal Y|_F$
isomorphe \`a $\mathcal X_F\times U$, et que dans cette d\'ecomposition,
la mesure
\begin{equation} \label{eq:mesureproduit}
\mu_{\mathcal Y,v}|_{\pi^{-1}(U)}
    = \mu_{\mathcal X,v} \otimes \mu_{\mathcal B,v}|_U . 
\end{equation}

Pour toute place $v$ de~$F$,
il r\'esulte du corollaire
au th\'eor\`eme~\ref{theo:picard}
la relation entre facteurs locaux
\begin{equation} \label{eq:factlocaux}
 L_v(s,\Pic({\mathcal Y}_{\bar F})) = L_v(s,\Pic({\mathcal X}_{\bar F}))
L_v(s,\Pic({\mathcal B}_{\bar F})). 
\end{equation}

Alors, les \'equations~\eqref{eq:mesureproduit} et~\eqref{eq:factlocaux}
impliquent
que la restriction de la mesure de Tamagawa de $\mathcal X(\A_F)$ \`a 
$\prod_v V(F_v)$ s'\'ecrit comme le produit
$$ \mu_{\mathcal Y,\Sigma}|_{\prod_v V(F_v)} 
     = \mu_{\mathcal X,\Sigma}
            \otimes \mu_{\mathcal B,\Sigma}|_{\prod_v U(F_v)}. $$
Or, si $\overline{U(F)}$ est l'adh\'erence de $U(F)$
dans le produit $\prod_v U(F_v)$,
l'adh\'erence de $V(F)$ dans $\prod_v V(F_v)$
s'identifie \`a $\overline{\mathcal X(F)}\times \overline{U(F)}$.
Int\'egrons $\mu_{\mathcal Y,\Sigma}$ sur $\overline{V(F)}$ ;
en utilisant le lemme~\ref{lemm:tamaouvert}, on obtient
$$ \int_{\overline{\mathcal Y(F)}} \mu_{\mathcal Y,\Sigma}
     = \int_{\overline{\mathcal X(F)}} \mu_{\mathcal X,\Sigma}
   \times \int_{\overline{\mathcal B(F)}} \mu_{\mathcal B,\Sigma}. $$

L'\'equation~\eqref{eq:factlocaux} implique aussi
que pour $\Re(s)>1$,
$$ L_\Sigma(s,\Pic({\mathcal Y}_{\bar F}))
    = L_\Sigma(s,\Pic({\mathcal X}_{\bar F}))
        L_\Sigma(s,\Pic({\mathcal B}_{\bar F})). $$
Par suite, l'ordre du p\^ole en $s=1$ pour la fonction $L_\Sigma$
de $\mathcal Y$ est la somme des ordes des p\^oles pour $\mathcal X$
et $\mathcal B$, et donc
$$ L_\Sigma^*(1,\Pic({\mathcal Y}_{\bar F}))
= L_\Sigma^*(1,\Pic({\mathcal X}_{\bar F}))
        L_\Sigma^*(1,\Pic({\mathcal B}_{\bar F})). $$
Le th\'eor\`eme est donc d\'emontr\'e.
\end{proof}

\subsection{Torseurs trivialisants}


Le paragraphe~\ref{subsec:torsion}
a montr\'e que le ph\'enom\`ene de torsion des m\'etriques
ad\'eliques intervient naturellement dans nos constructions.
Cependant, la hauteur tordue n'est facile \`a calculer que lorsqu'il existe
des sections propres pour l'action du groupe.
L'existence de sections canoniques permet comme on l'a vu 
de disposer d'une fonction hauteur sur les points ad\'eliques.

Les torseurs trivialisants que nous introduisons ici ont pour fonction
de fournir --- au prix d'un changement de vari\'et\'e --- d'une droite
canonique de sections.

Dans ce paragraphe, nous nous pla\c cons sur un corps $F$.
Supposons que 
$\Pic^G(\mathcal X)\simeq\Pic(\mathcal X)\times X^*(G)$
est un groupe de type fini.

Soit $H$ un groupe alg\'ebrique sur $F$, $\mathcal X_1\ra\mathcal X$
un $H$-torseur qui induise par fonctorialit\'e covariante des torseurs
un isomorphisme $X^*(H)\ra\Pic(\mathcal X)$.
On suppose de plus que $\mathcal X_1$ est muni d'une action de $G$
qui rel\`eve l'action de $G$ sur $\mathcal X$ et qui commute \`a l'action de
$H$.
On peut construire un tel $\mathcal X_1$ en fixant
$\lambda_1,\ldots,\lambda_h$ des fibr\'es inversibles $G$-lin\'earis\'es
sur $\mathcal X$ dont les classes forment une base de $\Pic(\mathcal X)$.
On pose alors $\mathcal X_1=\prod_{i=1}^h
(\V(\lambda_i^\vee)\setminus\{0\})$ et $H=\gm^h$.

Soit $T$ le plus grand quotient de $G$ tel que l'homomorphisme
naturel $X^*(T)\ra X^*(G)$
est un isomorphisme. (C'est le quotient de $G$ par l'intersection des 
noyaux des caract\`eres de $G$).
On pose $\tilde {\mathcal X}=\tilde {\mathcal X}_1\times T$
et $\pi:\tilde {\mathcal X}\ra \mathcal X$
la composition de la premi\`ere projection de de la projection $\tilde
{\mathcal X}_1\ra\mathcal X$. C'est un $H\times T$-torseur muni
d'une action de $G$ (diagonale).

\begin{enonce}[remark]{Exemple}
Supposons que $\mathcal X=P\backslash G$ est un espace de drapeaux g\'en\'eralis\'e
pour un groupe alg\'ebrique simplement connexe semi-simple $G$ sur $F$.
On a $\Pic(\mathcal X)\simeq X^*(P)$ et $G\ra \mathcal X$ est
un $P$-torseur qui induit un isomorphisme $X^* (P)\simeq \Pic(\mathcal X)$.
De plus, $T=\{1\}$. Ainsi, on peut prendre $\tilde {\mathcal X}=G$.
\end{enonce}

\begin{enonce}[remark]{Exemple}
Lorsque le groupe $G$ est trivial, on retrouve
les torseurs universels introduits dans le contexte des hauteurs
par Salberger et Peyre (cf.~\cite{salberger98}, \cite{peyre98}).
\end{enonce}

\begin{enonce}{Fait}
Si $\lambda\in\Pic^G(\mathcal X)$, $\pi^*\lambda$ admet une $F$-droite
canonique de sections $G$-invariantes.
\end{enonce}

\begin{rema}
L'isomorphisme canonique $\Pic^G(\mathcal X)\simeq X^*(H\times
T)=X^*(H)\times X^*(G)$ admet une r\'eciproque qu'il est facile d'expliciter.
En effet, soient $\chi_H$ et $\chi_G$ deux caract\`eres de $H$ et $G$
respectivement. On d\'efinit un fibr\'e inversible $G$-lin\'earis\'e sur $\mathcal
X$ comme suit : on quotiente
$\tilde {\mathcal X}\times\A^1 =\tilde {\mathcal X}_1\times T\times
\A^1$ par l'action de $H$ donn\'ee par
$$ h\cdot (\tilde x,t,u) = (h\cdot\tilde x,t,\chi_H(h)^{-1}u), \quad h\in H,
\quad (\tilde x,t,u)  \in \tilde{\mathcal X}_1\times T\times \A^1 $$
et la $G$-lin\'earisation provient de l'action de $G$ sur
$\tilde{\mathcal X}\times \A^1$ fournie par
$$ (\tilde x,t,u)\cdot g = (g\cdot \tilde x ,g\cdot t, \chi_G^{-1}(g)u),
 \quad g\in G,
 \quad (\tilde x,t,u)  \in \tilde{\mathcal X}_1\times T\times \A^1 . $$
\end{rema}

Par la construction~\ref{consY}, on obtient ainsi un $F$-sch\'ema
$\tilde {\mathcal Y}$ avec une projection $\tilde {\mathcal Y}\ra \mathcal Y$.
Supposons que $\mathcal Y$ provient de la situation~\ref{hatsitu},
on dispose en particulier de fibr\'es inversibles sur $\mathcal Y_ F$
munis de m\'etriques ad\'eliques associ\'es aux fibr\'es inversibles $(G,\mathbf
K)$-lin\'earis\'es sur $\mathcal X$. En particulier, on obtient sur $\tilde
{\mathcal  Y}$
des fibr\'es inversibles avec m\'etriques ad\'eliques. 
Le fait nouveau est que l'on dispose d'une hauteur sur
les points ad\'eliques de $\tilde{\mathcal Y}$ associ\'ee \`a ces fibr\'es inversibles.
En effet, une fois remont\'es \`a $\tilde{\mathcal Y}$,
ces fibr\'es inversibles poss\`edent une droite de sections
$F$-rationnelle canonique.

\subsection{Exemples}
\label{subsec:exemples}

\paragraph{Action d'un tore}

Pour les applications auxquelles notre deuxi\`eme article
sera consacr\'e, on consid\`ere l'action d'un tore $T$.

Un tel tore peut agir non seulement sur des vari\'et\'es toriques,
mais aussi sur des vari\'et\'es de drapeaux g\'en\'eralis\'ees $P\backslash G$,
via un morphisme $T\ra G$.

Dans le cas des vari\'et\'es toriques sur un corps de nombres $F$,
on dispose de mod\`eles canoniques
sur $\Spec\mathfrak o_F$ (si le tore est d\'eploy\'e),
et de m\'etriques hermitiennes
\`a l'infini canoniques sur les fibr\'es en droites.
Pour tout plongement $\sigma$ de $F$ dans $\C$, les points complexes $T(\C)$
du tore admettent un unique sous-groupe compact maximal $K_\sigma$,
et les m\'etriques hermitiennes introduites sont  automatiquement
$K_\sigma$-invariantes. On obtient ainsi des fibr\'es hermitiens
$(T,\mathbf K)$-lin\'earis\'es (cf.\  par exemple~\cite{batyrev-t95b}).

Dans le cas des vari\'et\'es de drapeaux $P\backslash G$,
une fois fix\'e des sous-groupes compacts maximaux de $G$ aux places \`a l'infini,
il est aussi possible de munir
les fibr\'es en droites $P$-lin\'earis\'es de m\'etriques hermitiennes
invariantes pour ces sous-groupes compacts maximaux
et donc pour le sous-groupe compact maximal de $T(\C)$.
Aux places finies, les m\'etriques $v$-adiques qu'on obtient
admettent une description analogue en termes de la d\'ecomposition
d'Iwasawa (cf.~\cite{franke-m-t89}).

D'autre part, un $T$-torseur sur un $F$-sch\'ema $\mathcal B$,
du moins quand le tore est d\'eploy\'e, \`a la donn\'ee d'un morphisme
$X^*(T)\ra\Pic(\mathcal B)$, et donc, une fois fix\'e une base de $X^*(T)$,
\`a des fibr\'es en droites $\lambda_1,\ldots,\lambda_t\in\Pic(\mathcal B)$.
(On a not\'e $t=\dim T$.)
La trivialisation des $T/K_\sigma$-torseurs
correspond, ainsi qu'on l'a dit apr\`es la d\'efinition~\ref{defi:torsarith} d'un
$T$-torseur arithm\'etique,
\`a une m\'etrique hermitienne sur les fibr\'es en droites $\lambda_i$.

Dans le cas o\`u $T$ agit sur une vari\'et\'e torique, on obtient
alors par la construction~\ref{hatsitu} une famille de vari\'et\'es
toriques sur $\mathcal B$.
On peut notamment compactifier ainsi une vari\'et\'e semi-ab\'elienne
$\mathcal T\ra\mathcal B$
et construire sur la compactification $\mathcal Y$ des fonctions hauteurs
canoniques. Dans ce cas, les $\lambda_i$ sont des fibr\'es en droites
alg\'ebriquement \'equivalent \`a $0$ sur une vari\'et\'e ab\'elienne $\mathcal B$.
Si on a pris soin de les munir, ainsi que tous les fibr\'es en droites
sur $\mathcal B$, de leur m\'etrique ad\'elique canonique,
pour laquelle le th\'eor\`eme du cube est une isom\'etrie,
on obtient sur $\mathcal Y$ les hauteurs canoniques,
au sens de la hauteur de N\'eron--Tate.
(Dans ce cas particulier, cf.~\cite{chambert-loir95} o\`u
l'on trouvera cette construction dans un esprit analogue,
et~\cite{cohen86}, o\`u est donn\'ee une construction {\og \`a la Tate\fg}
de ces hauteurs canoniques, due \`a M.~Waldschmidt).

Dans le cas o\`u $T$ agit sur une vari\'et\'e de drapeaux g\'en\'eralis\'ee,
on obtient la vari\'et\'e de drapeaux (g\'en\'eralis\'ee) d'un fibr\'e vectoriel 
sur $\mathcal B$ construit naturellement \`a partir des $\lambda_i$.
Ce cas \'etait \'etudi\'e (lorsque la base est aussi une vari\'et\'e de drapeaux)
dans la th\`ese de M.~Strauch (\cite{strauch97}).

\paragraph{Vari\'et\'es de drapeaux}

Tout fibr\'e vectoriel sur $\mathcal B$ donne lieu \`a des vari\'et\'es
de drapeaux g\'en\'eralis\'ees. Dans ce cas, le groupe $G$ est le groupe
lin\'eaire $\GL(d)$, $\mathcal X$ est une vari\'et\'e $P\backslash G$.
On identifie en effet un fibr\'e vectoriel de rang~$n$ sur $\mathcal B$
\`a un $\GL(d)$-torseur. Si l'on choisit comme sous-groupe compact
\`a l'infini le groupe unitaire $\mathrm U(d)$, 
la trivialisation \`a l'infini du
$G/K$-fibr\'e correspond \`a une m\'etrique hermitienne sur le fibr\'e vectoriel.

Il est \`a noter que cette situation se retrouve, mais dans l'autre sens,
dans le calcul du comportement de la fonction z\^eta des hauteurs
d'une puissance sym\'etrique d'une courbe $\mathcal C$ de genre $g\geq 2$.
Dans ce cas en effet, si $d > 2g-2$, $\Sym^d\mathcal C$
est un fibr\'e projectif au-dessus de la jacobienne de $\mathcal C$
associ\'e \`a un fibr\'e vectoriel de rang $d+1-g$.

\paragraph{Action d'un groupe vectoriel}

Dans~\cite{chambert-loir96} et~\cite{chambert-loir95}, on \'etudie
des compactifications d'extensions vectorielles de vari\'et\'es ab\'eliennes.
Expliquons comment ce travail s'ins\`ere dans les constructions de cet
article lorsque, pour simplifier les notations, on prend $G=\ga$.

Un $\ga$-torseur sur $\mathcal B$ correspond \`a une
extension de $ \mathcal O_{\mathcal B}$ par lui-m\^eme, soit un
fibr\'e vectoriel $\mathcal E$ de rang~$2$ sur $\mathcal B$.
La trivialisation du $\ga$-torseur \`a l'infini
correspond \`a un scindage $\mathcal C^\infty$ de l'extension
sur $\mathcal B(\C)$. D'autre part, $\ga$ agit naturellement sur $\P^1$
(via son plongement dans $\GL(2)$, $a\mapsto \left(\begin{smallmatrix} 1 &
a \\ 0 & 1 \end{smallmatrix}\right)$).
On obtient ainsi une compactification du $\ga$-torseur en une famille
de droites projectives sur $\mathcal B$.


\edef\guillemotleft{\guillemotleft\nobreak\,}
\edef\guillemotright{\unskip\nobreak\,\guillemotright}
\frenchspacing
\bibliographystyle{smfplain}
\bibliography{acl}

\end{document}